\documentclass[12pt]{scrartcl}

\usepackage[utf8]{inputenc}
\usepackage[T1]{fontenc}
\usepackage{placeins}
\usepackage{lmodern}
\usepackage{microtype}
\usepackage{graphicx}
\usepackage{caption}
\usepackage{amsmath}
\usepackage{amsfonts}
\usepackage{amssymb}
\usepackage{geometry}
\usepackage{url}
\usepackage{hyperref}
\usepackage{xspace}
\usepackage{color}
\usepackage{marginnote}
\usepackage{overpic}

\DeclareCaptionFormat{marglabel}{%
  \reversemarginpar
  \marginnote[\textbf{#1#2}]{\textbf{#1#2}}%
  #3%
}

\newcommand{\mesh}{$\mathcal{M}$\xspace}
\newcommand{\V}{\mathcal{V}}
\newcommand{\E}{\mathcal{E}}
\newcommand{\F}{\mathcal{F}}
\title{On The Topology of Polygonal Meshes}
\author{Andreas Bærentzen}
\date{\today \hspace{0.5cm}}

\begin{document}

\maketitle

\begin{abstract}
    This paper is an introductory and informal exposition on the topology of polygonal meshes. We begin with a broad overview of topological notions and discuss how homeomorphisms, homotopy, and homology can be used to characterise topology. We move on to define polygonal meshes and make a distinction between intrinsic topology and extrinsic topology which depends on the  space in which the mesh is immersed. A distinction is also made between quantitative topological properties and qualitative properties. Next, we outline proofs of the Euler and the Euler-Poincaré formulas. The Betti numbers are then defined in terms of the Euler-Poincaré formula and other mesh statistics rather than as cardinalities of the homology groups which allows us to avoid abstract algebra. Finally, we discuss how it is possible to cut a polygonal mesh such that it becomes a topological disc.
\end{abstract}

\section{Introduction}
\begin{figure}
    \centering
    \includegraphics[width=0.49\textwidth]{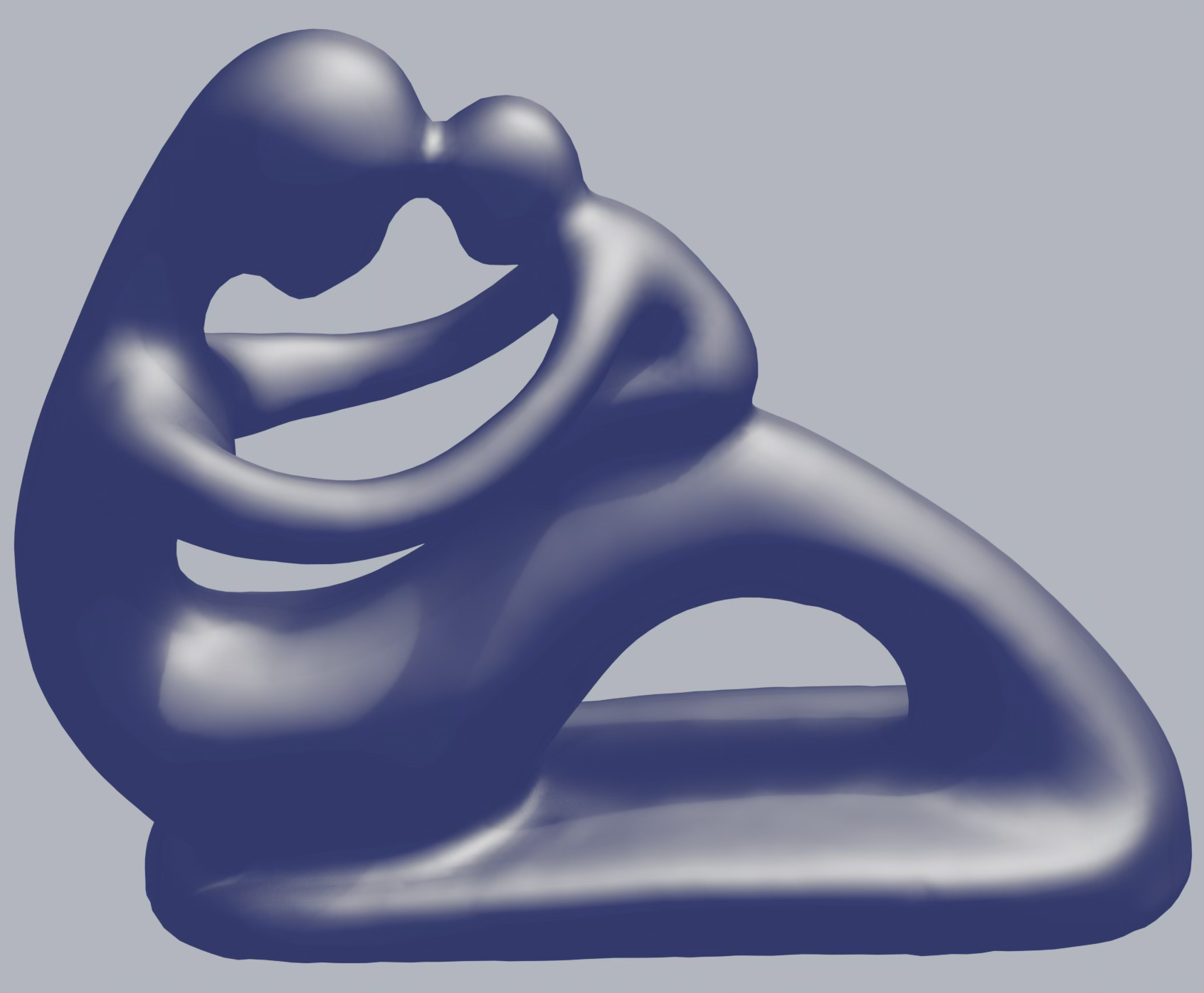}
    \includegraphics[width=0.49\textwidth]{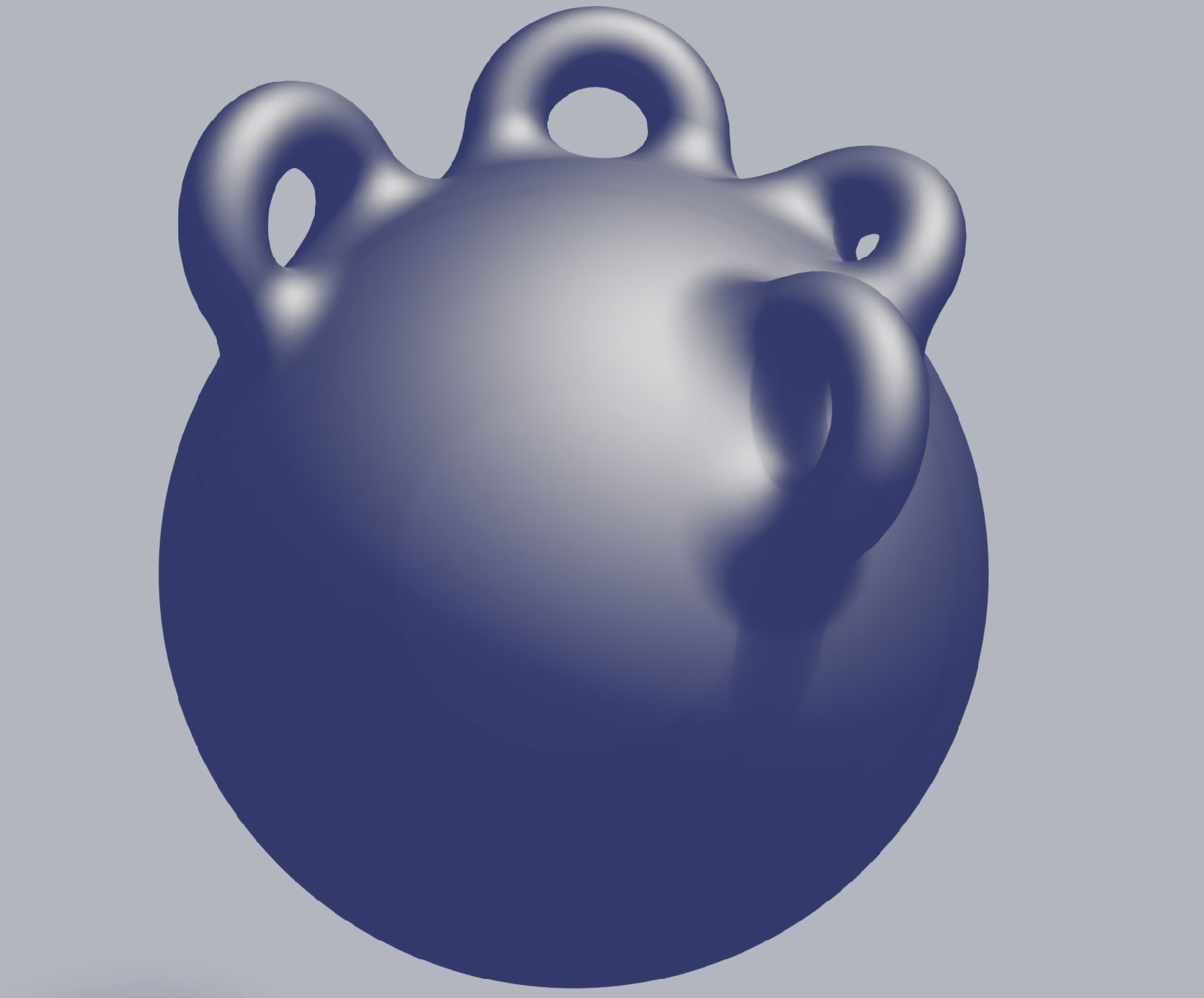}
    \caption{The Fertility model is shown on the left. It may not be immediately obvious how many holes (or handles) it possesses, but on the right we see a sphere with four handles which is topologically equivalent, i.e. there is a homeomorphism between these two shapes, and on the right it is obvious that the model has four handles.}
    \label{fig:fertility}
\end{figure}

Topology is the study of invariants of shapes under deformation. These invariants are quantities that do not change no matter how much you deform the shape -- as long as it is neither torn nor glued. A specific invariant that we often think about is the number of holes, and a classical example is that a cup with a single handle is topologically equivalent to a donut.

From a practical point of view, topological properties can be of great interest. For instance, we may want to count the number of objects in a 2D or 3D image; this is such a common task that it is rarely recognised as being related to topology. In other cases we may want to determine whether objects contain cycles or have a purely branching structure, or we may be interested in knowing whether an object has internal cavities.

In some situations we may be interested in topological properties because they have some biological significance \cite{Arnavaz2022}. In other situations we may need to compute topological properties in order to define loss functions for machine learning \cite{clough2020topological}. Finally, there are situations where we need to restrict the topology and want to detect potential topological changes in order to avoid them \cite{Cui2024}.

While there are many applications of topology, it can be a challenging topic unless one has a relatively strong mathematical background. Specifically, it is important to understand group theory \cite{alexandroff2012introduction}, which belongs to the field of abstract algebra, since the important Betti numbers are defined as cardinalities of homology groups. In this document, we sidestep some of the mathematical requirements by considering the topology of 3D surfaces only in the form of polygonal meshes and derive all quantitative topological properties from the Euler-Poincaré formula. By limiting the scope in this way, we arrive at formulas for key topological properties such as genus and Betti numbers without referring to abstract algebra. For readers who do want to understand homology, the book by Edelsbrunner and Harer \cite{edelsbrunner2010computational} is a good introduction that focuses on the computational aspects, and Giblin's book \cite{giblin2013graphs} is recommended as a more mathematical text. Readers interested in a deeper understanding of the Euler-Poincaré formula should consider Imre Lakatos' extraordinary book ``Proofs and Refutations'' \cite{lakatos1963proofs} which uses the formula merely as an example, but since it is the only example it is covered very thoroughly.
\subsection{Why is topology hard?}
Topological properties feel familiar and appeal to our intuition, but complications arise when we make the notions precise. Consider a plastic bag that has been sealed. Ignoring the slight thickness of the material, a plastic bag is a 2-dimensional surface which \textit{may} bound an interior. \marginnote{\footnotesize\itshape In the popular math book ``Shape'' \cite{ellenberg2022shape}, Ellenberg explains quite vividly just how confusing it can be to decide how many holes a straw has.} If we cut a hole in the bag with a pair of scissors it no longer bounds. In fact, it could be spread out to form a disc. Confusingly, the disc does not have a hole although, of course, its boundary is the original hole. In contrast, if we poke a hole through a lump of clay, we create a tunnel through the material that does not affect whether its surface bounds an interior. Effectively, we have then created a donut-like shape (i.e. a topological torus cf. Figure~\ref{fig:homs}) and the surface of this donut still bounds an interior. Thus, the informal term ``hole'' covers several distinct topological notions.

Another pernicious source of confusion is the fact that the same object can be understood in different ways, and these different ways may lead to different conclusions regarding its topology. This will be quite relevant in the following when we disentangle intrinsic topology which is indifferent to how the object is immersed in space and extrinsic topology which does take the spatial arrangement into account. 
To give another example of the importance of the framing, consider an egg-shaped object that is hollow inside. We might consider this object to consist of two surfaces: the outer surface and the surface that encloses the interior void. However, if we restrict our topological analysis to surfaces only, there is no way to distinguish an object with an interior cavity from two objects where one just happens to be inside the other. Since we will only be dealing with polygonal meshes in the following, it is worthwhile to be quite clear that polygonal meshes are surfaces, and there is no way of telling (from the representation itself) whether two disjoint polygonal meshes are part of the boundary of the same solid. 

\subsection{Background}
\begin{figure}
    \centering
    \begin{overpic}[width=\textwidth]{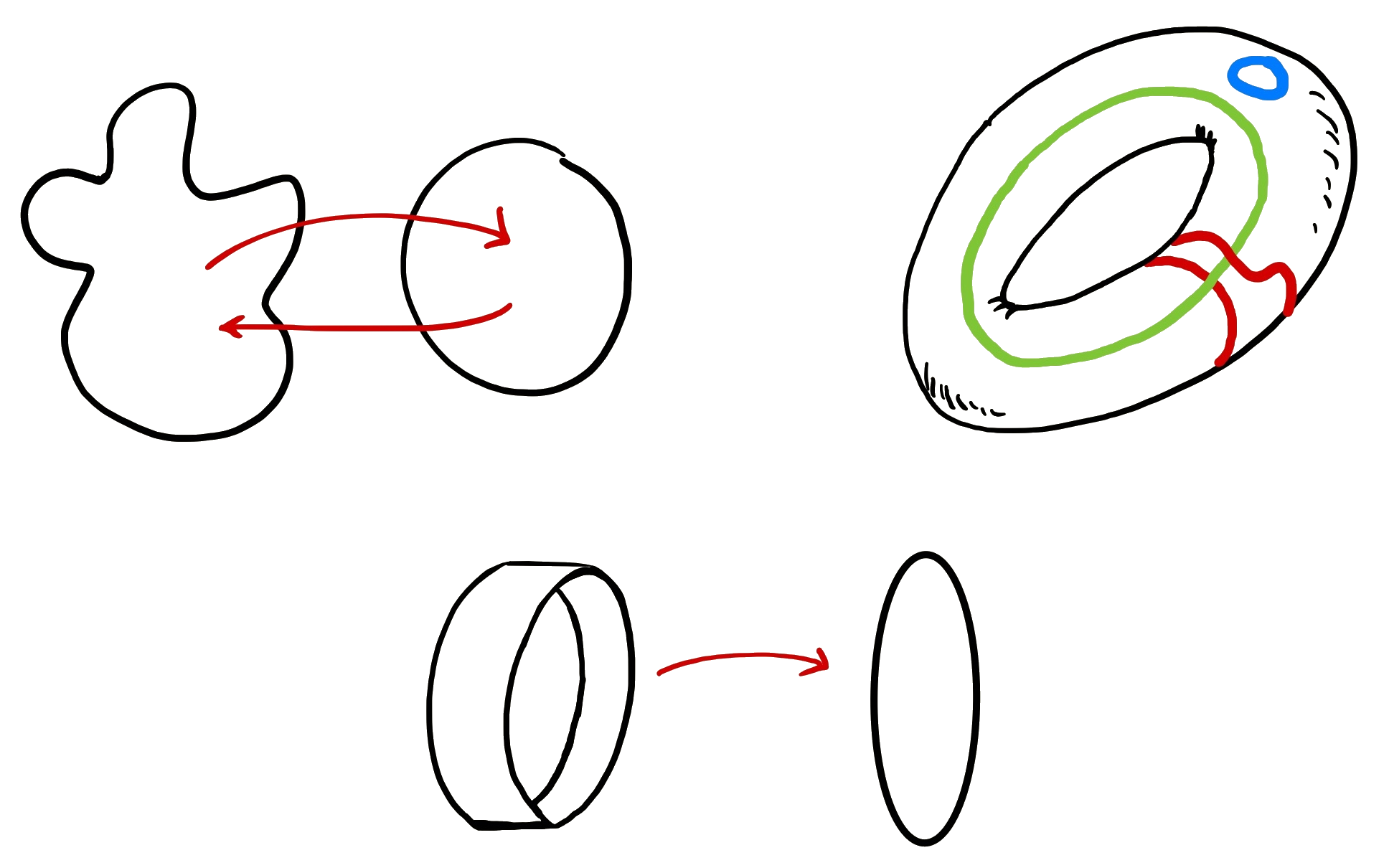}
        \put(25,49){$f$}
        \put(23,35){$f^{-1}$}
        \put(10,40){$A$}
        \put(35.5,41.5){$B$}
    \end{overpic}
    \caption{Top left: A homeomorphism, $f: A \rightarrow B$, is a continuous, bijective map which has a continuous inverse. Middle bottom: A homotopy is a continuous deformation. In this example an annulus is homotopic to a ring since the former can be continuously shrunk into the latter. Top right: The red curves are homologous because they differ only by the boundary of a region. However, the red curves are not homologous to the green or blue curves. The blue curve, on the other hand, is homologous to zero because it is the boundary of a region.}
    \label{fig:homs}
\end{figure}
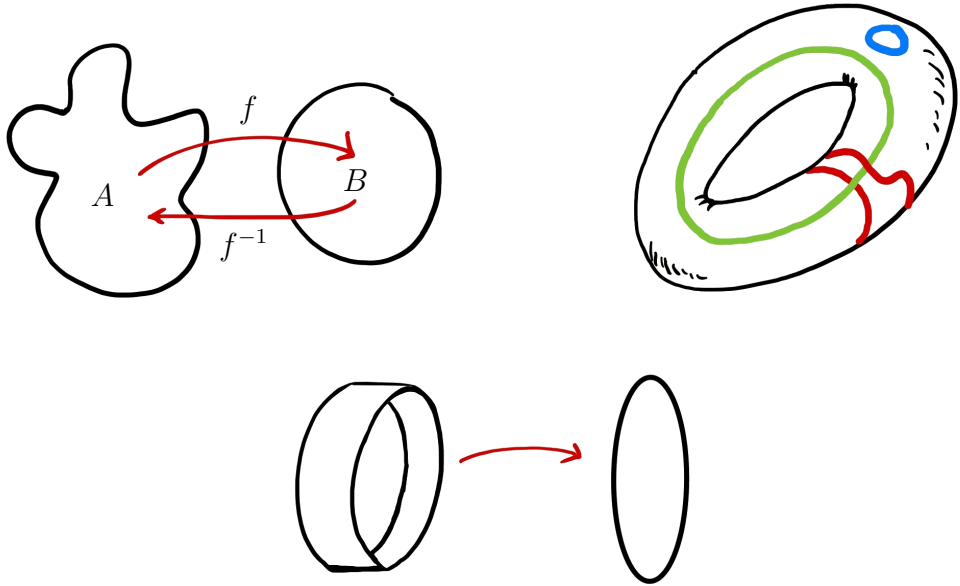
\marginnote{\footnotesize\itshape Note that we are here skipping over definitions from point-set topology. Also, we will only be concerned with $d=3$ in this text.} Before we turn to polygonal meshes, we will discuss some fundamental ideas that are key to understanding the topology of surfaces in general. First off, we need the notion of a \textit{topological space}. We will only consider topological spaces, $X\subset \mathbb R^d$, which are subsets of the normal Euclidean space and thus inherit the basic connectivity of $\mathbb R^d$. In the following, we will often refer to the topology of ``shapes'' or ``objects''. These two words are simply used as more familiar names for topological spaces.

\paragraph{Homeomorphisms} Two topological spaces are said to be \textit{homeomorphic} if there is a continuous bijective map from one to the other and this map has a continuous inverse (cf. Figure~\ref{fig:homs} top left). This is the strictest test of topological equivalence, but given two objects of identical topology there is no general method for computing a homeomorphism between them. Moreover, homeomorphisms only allow us to make statements about topological equivalence, they do not allow for quantification of topological properties. Nevertheless, we will often need to make statements about topological equivalence in terms of homeomorphisms in what follows. For instance, if we say that something is a topological disc, it means that there is a homeomorphism between the object and a disc.

\paragraph{Manifolds} 
\marginnote{\includegraphics[width=4.2cm]{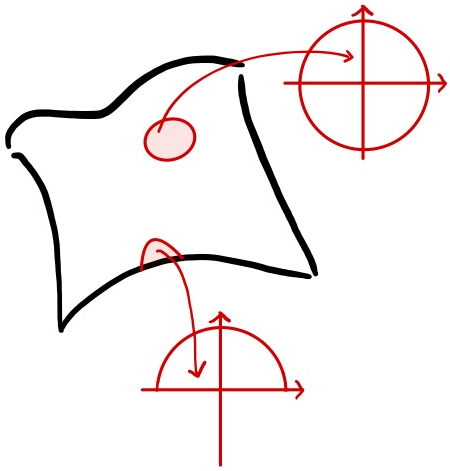}}
A topological space, $X$, is a \textit{2-manifold} if any sufficiently small neighbourhood of a point, $\mathbf p \in X$, is homeomorphic to an open disc (see the margin figure on the right). Intuitively, this means that we can cut out a (possibly tiny) patch which contains $\mathbf p$ and which can be flattened and stretched to a perfect disc. Thus, a 2-manifold models $\mathbb R^2$ locally at every point and can be seen as a mathematical definition of the notion of a closed surface. Of course, some surfaces, e.g. a sheet of paper, are open, and an open surface is a 2-manifold with boundary. Points on the boundary are homeomorphic to half-open disc as shown on the right.

The notion of a 2-manifold generalises beyond surfaces in 3D. An $n$-manifold in $\mathbb{R}^d$ where $n\le d$ is a set of points that is locally homeomorphic to either an open $n$-ball or a half-open $n$-ball in a sufficiently small region around every point belonging to the $n$-manifold. Because we will only consider surfaces, manifolds will always refer to 2-manifolds in the following.

\paragraph{Homotopy} Two maps, $f,g :X\rightarrow Y$, between topological spaces are considered \textit{homotopic} if one can be continuously transformed into the other, and the continuous transformation, $H: X \times [0,1] \rightarrow Y$ where $H(\mathbf x, 0)=f(\mathbf x)$ and $H(\mathbf x, 1) = g(\mathbf x)$, is called an \textit{homotopy}.

For our purposes it is enough to note that two topological spaces are considered \textit{homotopy equivalent} if one can be continuously deformed into the other (cf. Figure~\ref{fig:homs} middle bottom). At first glance, it may appear to be the same thing for two objects to be homeomorphic and homotopy equivalent, but that is actually false. A solid sphere is homotopy equivalent to a single point since we can smoothly contract the sphere to a point, but there is clearly no way to establish a bijective continuous mapping with a continuous inverse from a solid sphere to a single point. Conversely, we can cut a torus, tie it into a knot, and then stitch the ends back up. The knotted torus will still be homeomorphic to the original, but it is not homotopy equivalent since the continuous deformation would get stuck on the knot.

Shapes that belong to the same equivalence class are said to be of the same \textit{homotopy type}. For instance, this is the case for the sphere and the point or the annulus and the ring in Figure~\ref{fig:homs}.

\paragraph{Simplicial Homology} The notion of \textit{homology} from the field of algebraic topology provides us with a tool for quantifying topological properties. We will briefly outline the ideas of \textit{simplicial} homology which assumes that the object is a discrete structure which we can think of as a triangle mesh. Homology allows us to compute the so-called \textit{Betti numbers}, which are a mathematically precise way to count cycles (c.f. Section~\ref{sec:cycles}) according to the dimension of the cycle, are defined as the cardinalities of the \textit{homology groups}, $H_k$, and computed according to a deceptively simple formula,
\begin{equation}
    \beta_k = |H_k| = \left| \frac{Z_k}{B_k} \right| ,
    \label{eq:betti-via-homology}
\end{equation}
where $Z_k$ is the group of all k-dimensional cycles, and $B_k$ are the k-dimensional boundaries. For instance, $Z_1$ consists of all closed curves (i.e. 1D cycles) on the shape under consideration. Likewise, $B_1$ is the group of all boundaries of 2-dimensional regions of the shape. 

The involved entities are \textit{groups} \cite{alexandroff2012introduction}. A group is defined as a set of mathematical objects endowed with an operator that combines two objects, producing a new object which must again belong to the set, i.e. fulfilling the \textit{closure} property. Perhaps the most familiar example of a group is the set of integers with addition as the operation, but collections of mesh entities called \textit{chains} also constitute groups. A chain can be understood as a vector of mesh entities that can be either vertices, edges, or faces. The operator that combines two chains simply adds the entities modulo two. Hence, an entity is in the resulting chain if it belongs to either but not both of the chains combined by the operator.

Given an appropriate equivalence relation, groups can be partitioned into equivalence classes called \textit{quotient groups}. If we consider, $H_1$, the quotient group (formed by the division-like notation in Eq. \eqref{eq:betti-via-homology}) should be understood as a partitioning of all closed chains of edges (i.e. the members of $Z_1$) by the equivalence relation that they differ only by the boundary of a region (any element of $B_1$). In the case of an object that has sphere topology, all closed chains of edges are boundaries of some region. This means that they all differ from the empty cycle by the boundary of a region and we say that such chains are homologous to zero. On a torus there are two distinct classes of chains of edges that do not bound any regions. One of these classes contains chains that go around the tube and the other contains chains that go around the hole of the torus. See Figure~\ref{fig:homs} (top right) where the two red curves go around the tube and the green curve goes around the inner hole. The number of generators is called the rank, $|H_1|$, of the group. Thus, $|H_1|=2$ for the torus. 

Given two different chains that belong to the same class there is some patch on the surface that has their combination as its boundary, and we say that the two chains \textit{differ} by this boundary. For instance, in Figure~\ref{fig:homs} (top right) the two red curves bound a ribbon around the torus and differ by the boundary of this ribbon.

This brief discussion was meant to provide only the gist of homology groups. A much more thorough discussion is provided by Edelsbrunner and Harer \cite{edelsbrunner2010computational} (and many others) but is outside the scope of this text which is rather to demonstrate how the Betti numbers can be computed \textit{without} making any reference to homology. Instead, we will demonstrate that an equivalent formulation of the Betti numbers can be derived from the Euler-Poincaré formula.
\section{Polygonal Meshes}
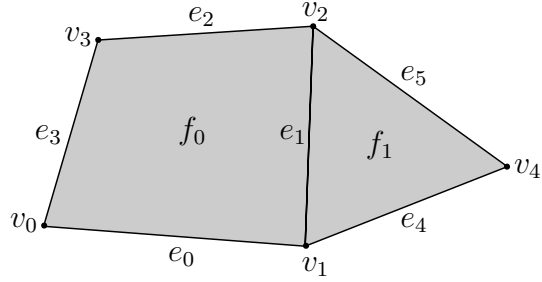
\begin{figure}
    \centering
    \begin{overpic}[width=0.5\textwidth]{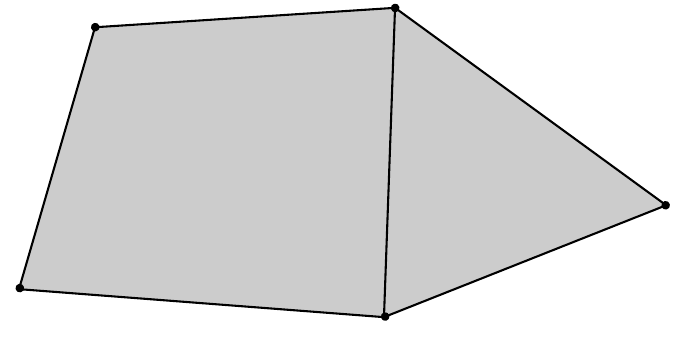}
        \put(-4,7){$v_0$}
        \put(55,-2){$v_1$}
        \put(55,50){$v_2$}
        \put(7,45){$v_3$}
        \put(98,18){$v_4$}
        \put(28,0){$e_0$}
        \put(50.5,25){$e_1$}
        \put(32,49){$e_2$}
        \put(1,25){$e_3$}
        \put(75,7){$e_4$}
        \put(75,37){$e_5$}
        \put(30,25){$f_0$}
        \put(68,22){$f_1$}
    \end{overpic}
    \caption{A polygonal mesh, $\mathcal M = \langle \V, \E, \F \rangle$, where $\V = \{v_0, v_1, v_2, v_3, v_4, v_5\}$, $\E=\{e_0, e_1, e_2, e_3, e_4, e_5\}$, and $\F = \{f_0, f_1\}$. $V=5$, $E=6$. and $F=2$. The boundaries of the two faces are $\partial f_0=\{e_0,e_1,e_2,e_3\}$ and $\partial f_1=\{e_1,e_4,e_5\}$. The edge boundaries are $\partial e_0=\{v_0,v_1\}$, $\partial e_1=\{v_1,v_2\}$, $\partial e_2=\{v_2,v_3\}$, $\partial e_3=\{v_0,v_3\}$, $\partial e_4=\{v_1,v_4\}$, and $\partial e_5=\{v_4,v_2\}$.}
    \label{fig:polygonal-mesh}
\end{figure}
\marginnote{\footnotesize\itshape By this definition, vertices that belong to an edge must be in the set of vertices and the edges of a face must belong to the set of edges. On the other hand, we may allow ``dangling'' edges or vertices. However, such unused entities are often considered defects.} 
A polygonal mesh can be understood as a generalisation of the (ancient) notion of a polyhedron. Like a polyhedron, a polygonal mesh is made of polygons that are joined along the edges. However, unlike a polyhedron, a mesh need not be closed, and the polygons are not necessarily perfectly flat.
Formally, a polygonal mesh is a tuple, \mesh = $\langle \V, \E, \F \rangle$, where $\V$ is a set of \textit{vertices}, $\E$ is a set of \textit{edges}, and $\F$ is a set of polygonal \textit{faces}. Vertices, edges, and faces are all considered \textit{mesh elements}. We will often need to refer to the number of vertices, $V=|\V|$, the number of edges, $E=|\E|$, and the number of faces, $F=|\F|$. We can think of vertices as 0-dimensional points in 3D and of edges as 1-dimensional line segments bounded by two vertices. Finally, faces are 2-dimensional polygons bounded by edges. A simple polygonal mesh consisting of just two {faces} is shown in Figure~\ref{fig:polygonal-mesh}. 
\subsection{Boundaries}
From a higher dimensional mesh element (e.g. a face) we can obtain the lower dimensional elements that bound it (e.g. a set of edges) by exploiting the \textit{boundary relationship}. The boundary relationship is formalised using the operator, $\partial$, which maps a set of $n$-dimensional mesh elements to their boundary which is a set of $(n-1)$-dimensional mesh elements. The boundary of a face is a set of edges,
\begin{equation}
\partial f = \{ e_0, ..., e_n \} \,,
\end{equation}
where $e_i \in \E$ and $f \in \F$. Likewise, the boundary of an edge is the set containing its two vertices,
\begin{equation}
\partial e = \{ u, v \} \,,
\end{equation}
where $u,v \in \V$ and $e \in \E$. Vertices do not have a boundary, hence
\begin{equation}
\partial v = \emptyset \,,
\end{equation}
where $v\in \V$. 

It is possible to talk about the boundary of a set of mesh elements, but first we need to define the (set-valued) co-boundary relation,
\begin{equation}
\delta_{\mathcal X}(\epsilon) = \{ x_i \in {\mathcal X} | \epsilon \in \partial x_i \} \,.
\end{equation}
\marginnote{\footnotesize\itshape The elements of $\mathcal X$ must be all of the same type -- either vertices, edges, or faces.}
For a set, ${\mathcal X}$, of mesh elements, $\delta_{\mathcal X}(\epsilon)$ is the subset of elements in ${\mathcal X}$ whose boundaries contain $\epsilon$. Now, we define the boundary of ${\mathcal X}$:
\begin{equation}
\partial {\mathcal X} = \{ \epsilon \in \bigcup_{x_i \in {\mathcal X}} \partial x_i | \left|\delta_{\mathcal X}(\epsilon)\right| = 1 \}
\end{equation}
In the case of faces, the boundary is thus the edges that belong to the boundary of precisely one face from the set. The analogous statement is true for a set of edges while for a set of vertices the boundary is empty.

Continuing the example from Figure~\ref{fig:polygonal-mesh}, if we let $\F= \{f_0, f_1\}$ denote both faces (i.e. the entire mesh), it is clear that $\partial \F = \{e_0, e_2, e_3, e_4, e_5\}$, since the shared edge $e_1$ is excluded from the boundary of $\F$.
\subsection{Cycles}
\label{sec:cycles}
\begin{figure}
    \includegraphics[width=\textwidth]{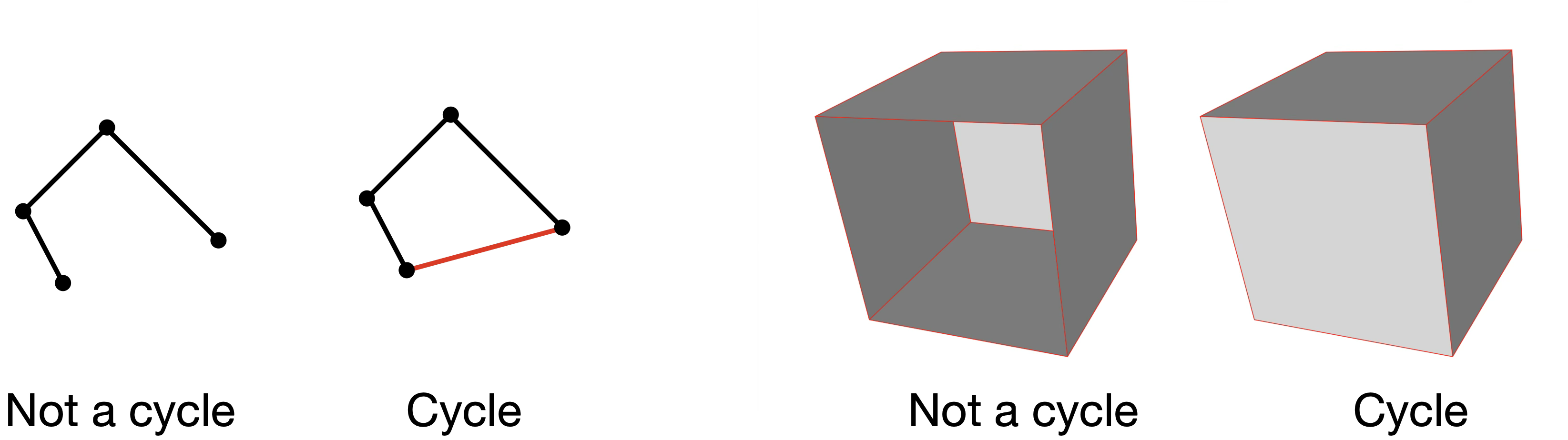}
    \caption{Left: two sets of edges where the left one is not a cycle, but the right one is a cycle. The cycle instigator is shown as a red line segment. Right: a cube which is not a cycle since it is missing a face is shown on the left. On the right a closed cube which constitutes a cycle. Both of the shown cycles are also simple.}
    \label{fig:cycle-examples}
\end{figure}
A \textit{cycle} is simply a set, $\mathcal C$, of mesh elements of the same dimension which has an empty boundary. In other words, if $\mathcal C$ is a cycle it is true that
\begin{equation}
\partial \mathcal C = \emptyset \,.
\end{equation}
Sets of faces or edges are \textit{closed} if they are cycles.

The notion of {cycles} is important in topology, and we can define cycles for sets of vertices, edges, and faces, although edge cycles are the most interesting. In fact, vertices and sets of vertices are trivially cycles.

A cycle, $\mathcal C$, is \textit{simple} if it consists of a single connected component and fulfils the following predicate:
\begin{equation}
\text{SimpleCycle}(\mathcal C) \iff
\forall \epsilon \in \bigcup_{x_i \in \mathcal C} \partial x_i : \left|\delta_{\mathcal C}(\epsilon)\right| = 2 \,.
\end{equation}
This means that a cycle of edges is simple if it is a connected set of edges and all vertices that are incident on the cycle are incident on precisely two edges. Note also that vertices and sets of vertices are not simple cycles.

The boundaries of faces of a polygonal mesh must be simple cycles:
\begin{equation}
\forall f \in \F: \text{SimpleCycle}(\partial f) \,.
\end{equation}
This is an important requirement and it is a matter of definition although, in a few cases, we may want to deviate slightly from the requirement that face boundaries must be simple (cf. Section~\ref{sec:computer-representation}).
\subsection{Manifoldness}
\label{sec:manifold}
Any point on a manifold has a small neighbourhood that is homeomorphic to a disc. This is trivially true for points inside faces since faces have disc topology by construction. To test a mesh, \mesh, for intrinsic manifoldness, we need to check that:
\begin{enumerate}
    \item Every edge of \mesh is adjacent to either one or two faces, i.e. \marginnote{\footnotesize\itshape If \mesh has no boundary curves then every face is adjacent to exactly two faces. In other words, a connected, closed, manifold mesh, \mesh, is a simple cycle.}
\begin{equation}
\forall e \in \E : \left| \delta_{\F}(e) \right| \in \{1,2\}
\end{equation}
    \item The mesh elements incident on a vertex, $v$, form only a single connected component, i.e. it must be true that
\begin{equation}
\forall e_0, e_1 \in \delta_{\E}(v) : e_0 \sim e_1
\end{equation}
where the transitive relation $\sim$ is defined by the conditions that 
\begin{eqnarray}
e_0 \sim e_1 &\impliedby& \delta_\F(e_0) \cap \delta_\F(e_1) \ne \emptyset \\
e_0 \sim e_1 &\impliedby& \exists e_2 \in \delta_{\E}(v) : e_0 \sim e_2 \wedge e_2 \sim e_1     
\end{eqnarray}
\end{enumerate}
If \mesh passes both of these tests, we know that it can everywhere be flattened to a disc (or half disc) and is thus manifold (possibly with boundary).

Provided a mesh is manifold we can prove that boundary edges always form cycles. As stated above, the faces incident on a vertex must form a single component. This means that a vertex is always incident on either two or zero boundary edges since the incident faces must form either a fan (as shown right) or a disc around the vertex. \marginnote{{\centering\includegraphics[width=3cm]{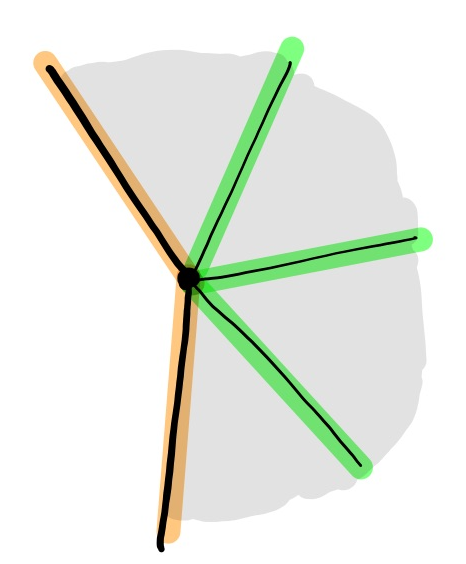}}\\ \footnotesize\itshape The faces around the center vertex form a fan capped by two (orange) boundary edges.} Consequently, all vertices which belong to boundaries are incident on two boundary edges (see the figure on the right). This means that we can always jump from a boundary edge to the next edge across one of its vertices, but, assuming the mesh is finite, we must eventually revisit the edge we started from. 

\begin{figure}
    \centering
    \includegraphics[width=0.5\textwidth]{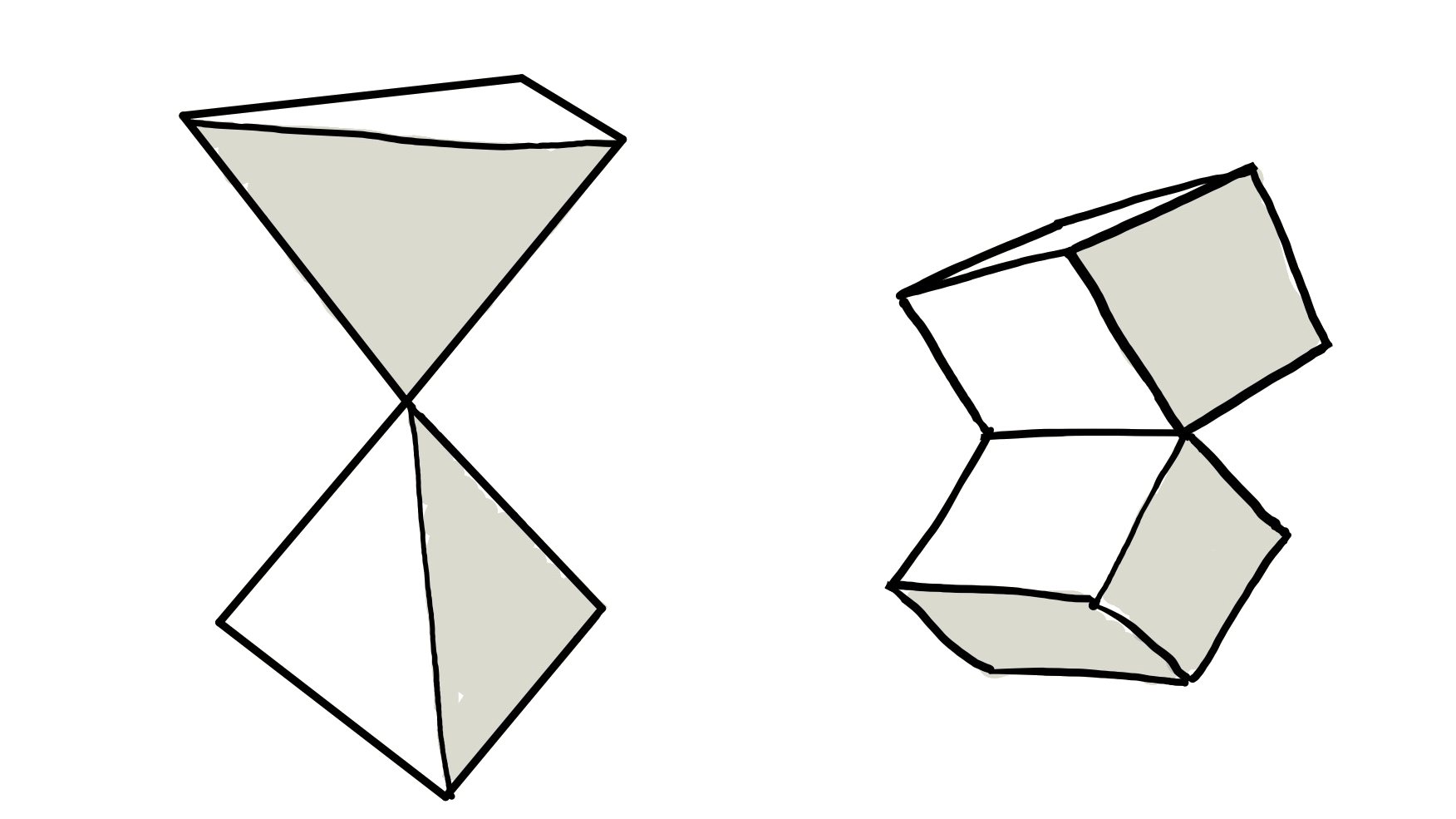}
    \caption{Two non-manifold configurations: on the left the two tetrahedra share a vertex, and on the right, two boxes share an edge. From the images we cannot tell whether the meshes are intrinsically or extrinsically non-manifold.}
    \label{fig:non-manifold}
\end{figure}

A mesh which is produced by joining two pyramids (one upside down) such that the tips of the pyramids are a shared vertex is not manifold. We could try to flatten a region around the shared vertex, but it would always end up as two discs connected at a single point. Likewise, if two cubes are placed side by side such that two edges touch (forming a shared edge) that is also a non-manifold configuration since four faces then meet at the edge where the cubes are joined. The two non-manifold configurations are shown in Figure~\ref{fig:non-manifold}.



\subsection{Computer Representation}
\label{sec:computer-representation}
From the perspective of computer representation, meshes can be stored in a number of ways that differ mostly in which of the three primitive types that is considered to be the primary one. 

Meshes can be stored as lists of faces where each face is represented as a cyclic list of vertices; two consecutive vertices then define an edge. This representation is preferred in computer graphics for rendering, and its simplicity makes it appealing also for disc storage of meshes. However, for manipulation of meshes, an edge-based data structure is often convenient since this type of representation makes certain editing operations efficient. For most edge-based data structures the requirement that an edge can be incident on either one or two faces is a built-in constraint. 

In the present context, the representation generally does not matter, but there is one further concern that should be addressed. While at least one popular data format stores the vertices directly as 3D points in space this is not a good idea if we are interested in topology. Vertices should be identified by indices, allowing us to easily find out if the corners of two polygons correspond to the same vertex (simply by checking whether the vertex indices are the same).

Of course, if we are interested in the geometry of the mesh, we need a separate table that maps vertex indices to points in space.
\paragraph{Acceptable Non-manifold Configurations}
One type of non-manifoldness is generally accepted in practice. 
Sometimes, a vertex is incident on more than two boundary edges. 
\marginnote{\centering\includegraphics[width=3cm]{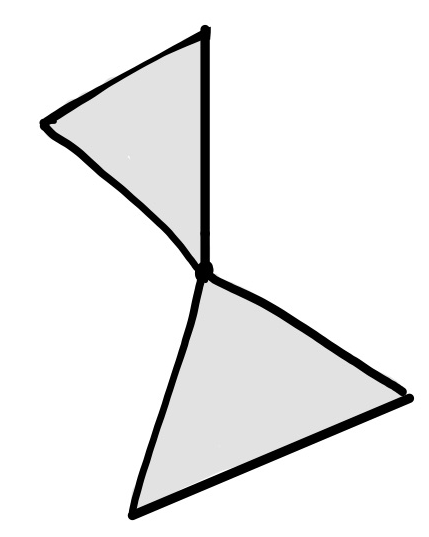}} 
An example is shown in the figure on the right where two triangles share a vertex. The result is that the vertex is incident on four boundary edges. As long as the mesh representation keeps track of which boundary edges that are adjacent in the incident boundary cycles, it is not a problem to have any even number of boundary edges incident on a vertex, but it is still required that the faces incident on the vertex all come in connected components that contribute two boundary edges. Put differently, the configuration is considered valid if it can be made manifold by closing the holes. 
\paragraph{Multiple Incidence} 
\marginnote{\includegraphics[width=4.2cm]{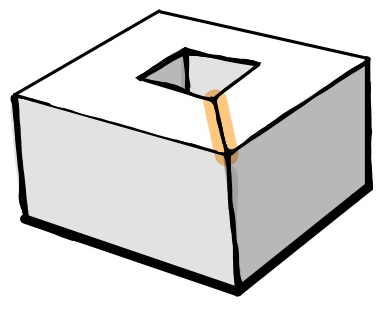} \footnotesize\itshape The highlighted edge is incident twice on the same face. Removing it would disconnect the interior boundary cycle from the exterior.} A special case that bears considering is that of \textit{multiple incidence}. An edge might be incident on the same face twice. Likewise, a vertex can be incident on the same face multiple times. It is probably clear that this breaks the simple-cycle requirement, but in practice mesh libraries can often handle such configurations.

Although multiple incidence is not necessarily problematic in itself, it can lead to problems if we manipulate a mesh in a computer. For instance, we might merge the two faces that are incident on a given edge. If the two faces are really the same as in the example shown on the right, we create two disconnected cycles that bound the same face. As discussed above, this should be avoided.
\paragraph{Triangle Meshes}
Triangle meshes are, of course, a special type of polygonal meshes where all faces have three edges. Numerous algorithms exist that operate only on triangle meshes and not general polygonal meshes. This is not a big limitation in practice since it is always possible to triangulate the faces of a polygonal mesh, thereby obtaining a triangle mesh. Moreover, data structures for triangle meshes are simpler than those for generic polygonal meshes.
\subsection{Intrinsic and Extrinsic Topology}
We have to face the fact that the topology of a polygonal mesh can mean two different things. 
\begin{itemize}
    \item If we are interested in the \textit{intrinsic} topology, the geometry is not important; a mesh is simply a schematic diagram of the structure of a surface. Since this diagram could be realised in 3D its topology is well defined, and this is what we call the intrinsic topology. 
    \item If we study the \textit{extrinsic} topology, the mesh is \textit{immersed} in 3D and this immersion \textit{is} the polygonal mesh. According to this view, a polygonal mesh is simply a set of points in 3D, namely the image of the mesh under the immersion, and the extrinsic topology is that of the image. Put differently, the extrinsic topology is that which is \textit{induced} by the immersion.
\end{itemize}

It is important to be aware of this distinction since the extrinsic and the intrinsic topology do not have to agree. Specifically, \mesh may intersect itself in the immersion which would often cause its extrinsic topology to differ from the intrinsic. For instance, a mesh may be intrinsically manifold but two vertices or two edges may touch as shown in Figure~\ref{fig:non-manifold}. In this case, the mesh is extrinsically non-manifold despite being manifold from the intrinsic point of view.

Some mesh configurations are valid if we are only interested in the intrinsic topology of abstract meshes, but make no sense when we consider a polygonal mesh with vertex positions in 3D and straight edges. For instance, from an intrinsic point of view, it makes sense that a vertex is connected to itself, but if we are interested in the extrinsic topology it is problematic since such an edge has no geometric meaning. Two vertices can also be connected by more than one edge. Extrinsically, this often makes the mesh non-manifold, but intrinsically it does not.
\subsection{Qualitative Topology}
\label{sec:qualitative}
The topological invariants \textit{watertightness} and \textit{orientability} are often important in the context of polygonal meshes. We can tell that they are topological properties since they are unaffected by a deformation of the mesh, but unlike the number of holes they are either true of false for a given mesh and thus qualitative rather than quantitative. 
\subsubsection{Watertightness}
\mesh is \textit{watertight} if it divides space into an exterior region, an interior region, and \mesh itself. If this is true, any path from a point in the interior to a point in the exterior contains at least one point that belongs to \mesh. A cube (and polyhedra in general) is watertight, but remove one of its faces (as shown in the figure on the right), and the property is lost. \marginnote{\centering \includegraphics[width=3cm]{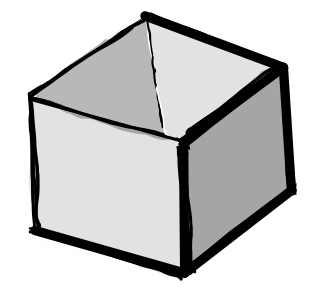} \footnotesize\itshape \\ Not watertight.} Again, this is the extrinsic definition of watertightness. Intrinsically, \mesh is watertight if it is closed, i.e. it has no boundary. A mesh which is not watertight has at least one cycle (cf. Section~\ref{sec:cycles}) of boundary edges.
\subsubsection{Orientability}
\mesh is \textit{orientable} if we can assign consistent orientations to all its faces. Even though the faces are not necessarily flat, a face does have a top side and a bottom side, but which is which is completely arbitrary. \mesh is orientable if we can assign top sides consistently.

How does one check whether the assignments are consistent? The trick is to use the induced circulation direction on the edges of the face. Assume we have assigned a top side to a face and have defined an arrow that sticks out of the face with the correct top-side orientation. We can now impose an orientation on the edges and vertices of the face (as shown in the figure on the right). Usually, we orient the edges counter clockwise with respect to the top-side direction -- consistent with the right hand rule: stick your right thumb in the top-side direction, and the other fingers will curl in the circulation direction. \marginnote{\includegraphics[width=4.2cm]{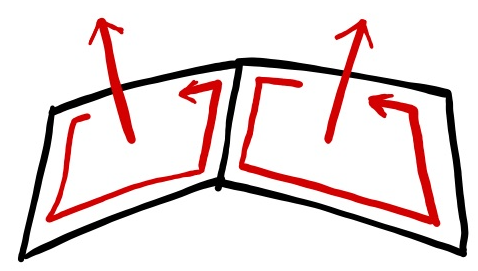}} \mesh is consistently oriented if every edge that is adjacent to two polygonal faces has opposite circulation direction with respect to these two faces (as shown in the figure on the right). Finally, \mesh is orientable if a consistent orientation \textit{can} be assigned. The procedure for testing orientability is to assign an arbitrary orientation to any vertex of \mesh and then recursively assign consistent orientations to its neighbours. \marginnote{\includegraphics[width=4.2cm]{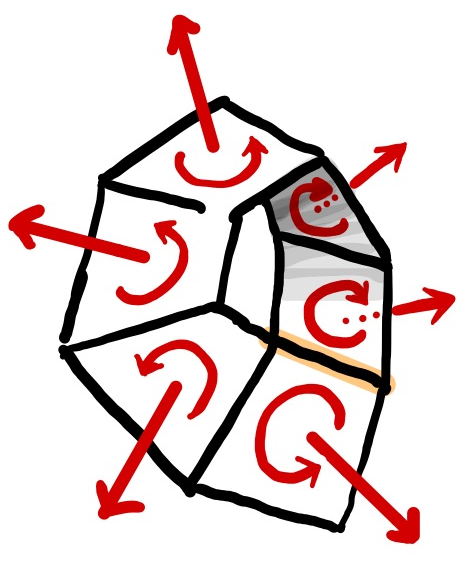}
\footnotesize\itshape The highlighted edge has the same direction with respect to both faces. On a Möbius band this problem will always arise.} If at some point this fails because two adjacent faces have been assigned inconsistent orientations, we know that the mesh is not orientable. An example is shown in the figure on the right.  

Our procedure for setting the orientation of a face is extrinsic because it depends on choosing a direction in the space into which it is immersed. However, the cyclic ordering of the edges of a face is an intrinsic property, and whether a consistent ordering of all faces can be achieved is likewise intrinsic.

There are several famous examples of non-orientable surfaces such as the Möbius band and the Klein bottle, but in the following we will generally assume that the polygonal meshes we consider are orientable. Non-orientable surfaces are confusing in the context of the topology of polygonal meshes. It is relatively clear that a Möbius band is not homeomorphic to an annulus, so our strongest tool can distinguish this non-orientable surface from its closest orientable relative, but they have the same homotopy type, and the quantitative methods that we will discuss later cannot make this distinction: we get the same numbers for an annulus and a Möbius band. On the other hand, it is easy to detect that a surface is non-orientable. We simply have to assign consistent orientations to all its faces and realise that this is impossible. In fact, some polygonal mesh representations (in the computer science sense) cannot represent non-orientable meshes, and a seam will be introduced somewhere. 
\section{Euler's Formula}
\begin{figure}
    \includegraphics[width=\textwidth]{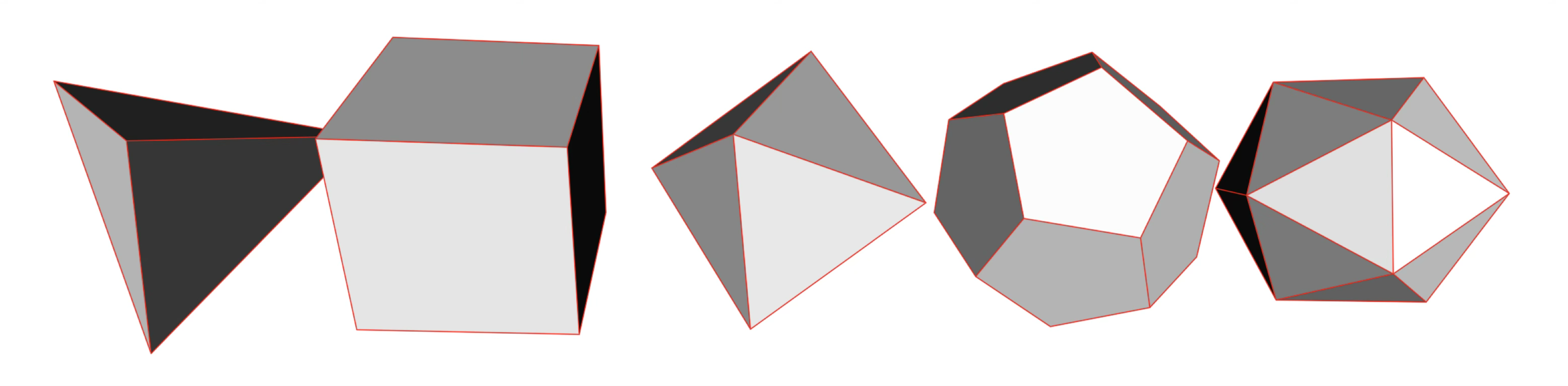}
    \caption{The Platonic solids are shown above. From left to right: the tetrahedron, cube, octahedron, dodecahedron, and icosahedron.}
    \label{fig:platonics}
\end{figure}

Our investigation into how we can quantify topological properties of polygonal meshes begins with Euler's formula,
\begin{equation}
    V-E+F=2\,,
    \label{eq:Euler}
\end{equation}
which holds for any polygonal mesh that is a topological sphere. \marginnote{\footnotesize\itshape \includegraphics[width=4.2cm]{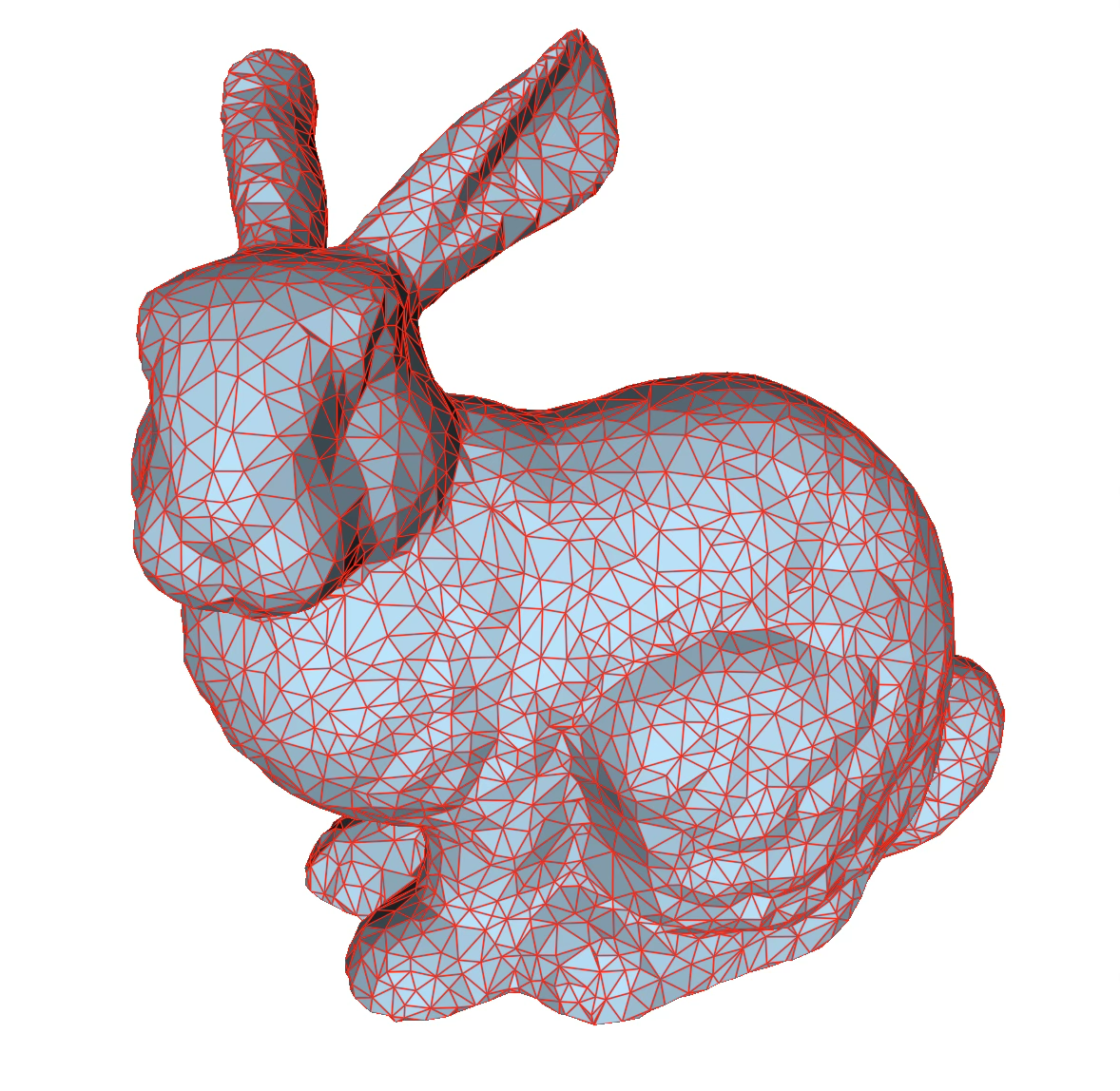} This mesh is a topological sphere.} A polygonal mesh is considered to be a topological sphere if it is homeomorphic to a sphere. Some of the most famous examples of polygonal meshes which are also topological spheres are the Platonic solids which are shown in Figure~\ref{fig:platonics}. By consulting Table~\ref{tab:example-meshes} we can ascertain that they all fulfil Eq.~\ref{eq:Euler}.
Of course, the Platonic solids are special because they have perfect symmetry in the sense that all faces and all vertices look exactly the same. This is not true of the mesh shown on the right. Yet the version of the Stanford Bunny that is shown in the margin is also a topological sphere, and if we plug the numbers from Table~\ref{tab:example-meshes} into Eq. \eqref{eq:Euler} it holds.

\begin{table}[ht]
\caption{Vertex ($V$), edge ($E$) and face ($F$) counts for the Platonic solids and the Stanford bunny. All examples have Euler characteristic $V-E+F=2$.}
\centering
\begin{tabular}{lrrr}
\hline
Name & $V$ & $E$ & $F$ \\
\hline
Tetrahedron & 4 & 6 & 4 \\
Cube & 8 & 12 & 6 \\
Octahedron & 6 & 12 & 8 \\
Icosahedron & 20 & 30 & 12 \\
Dodecahedron & 12 & 30 & 20 \\
Stanford Bunny & 3485 & 10449 & 6966 \\
\hline
\end{tabular}
\label{tab:example-meshes}
\end{table}

\begin{figure}
    \centering
    \includegraphics[width=0.49\textwidth]{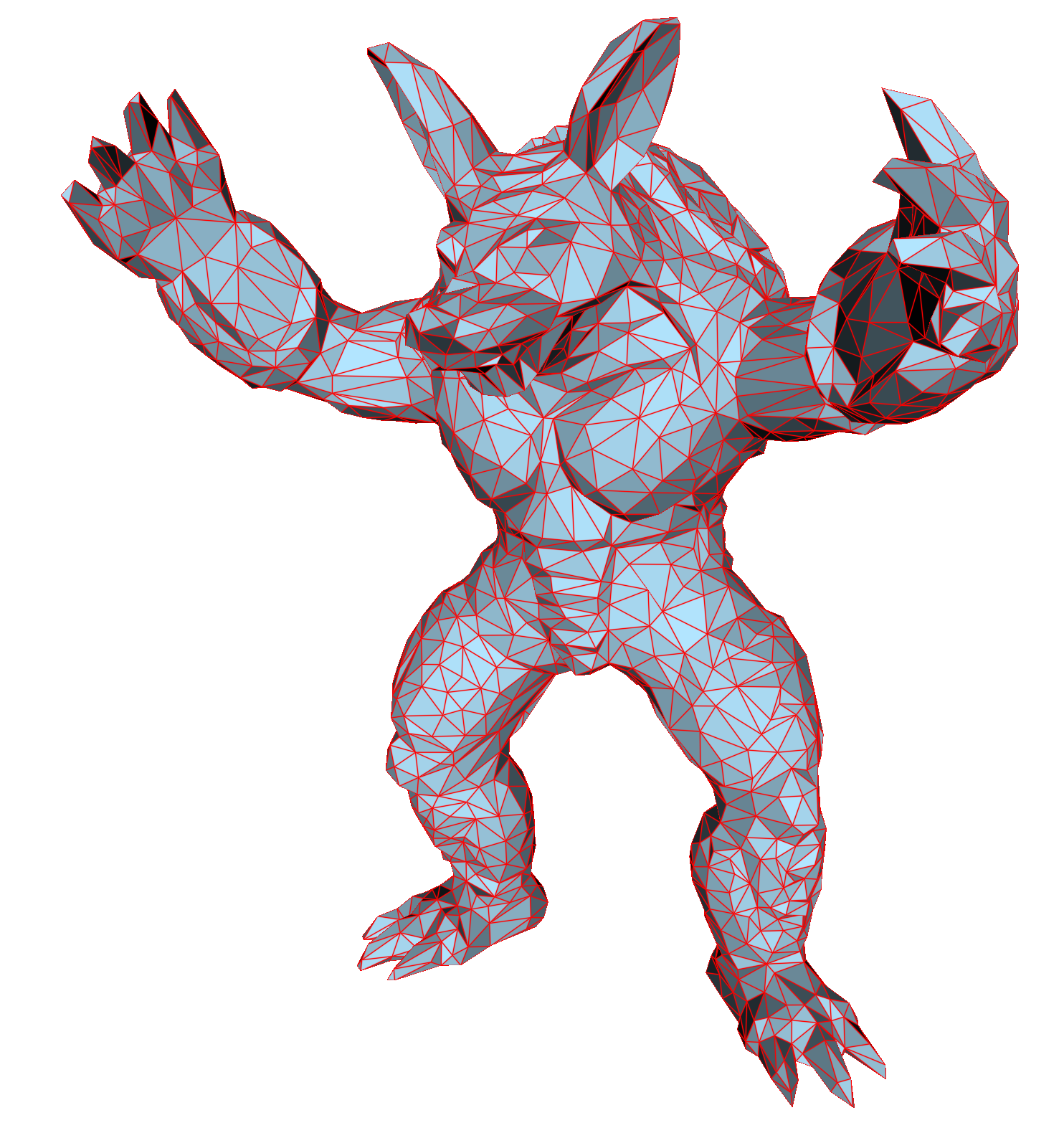}
    \includegraphics[width=0.49\textwidth]{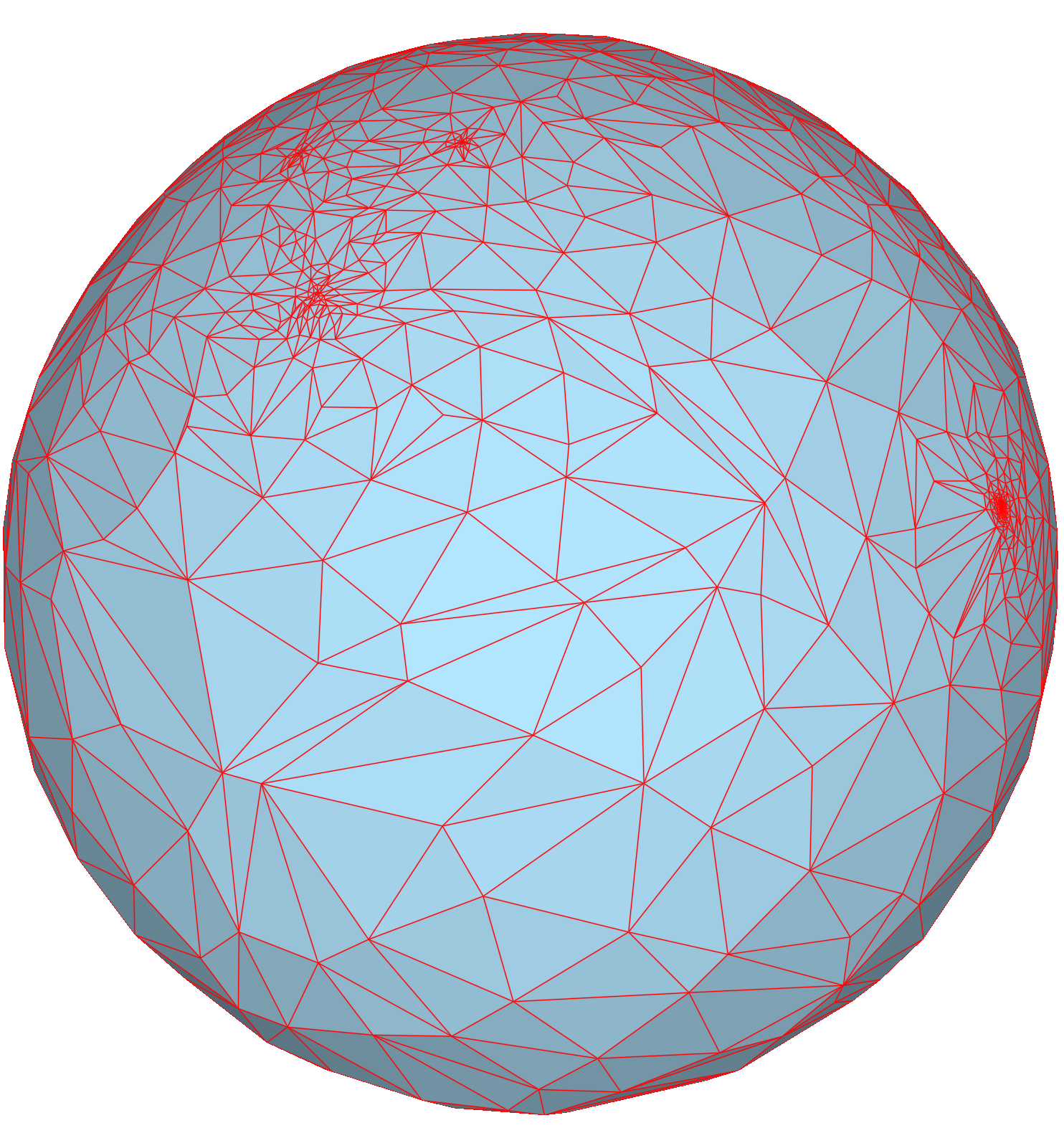}
    \caption{A simplified version of the Stanford Armadillo is shown left. On the right is a mesh consisting of the exact same vertices and faces but transformed such that all vertices lie on the sphere. The edges remain straight lines, but it is easy to imagine that they could be replaced with great arcs: if we imagine a light placed in the centre of the sphere, the shadows that the edges would cast on the (inside of) the sphere would be precisely the needed circular arcs as explained by Eppstein in his Geometry Junkyard.}
    \label{fig:armadillo-sphere}
\end{figure}

Euler's formula is arguably surprising, but it is fortunately easy to prove\marginnote{\footnotesize\itshape In fact, David Eppstein has 21 proofs in his
\href{https://ics.uci.edu/~eppstein/junkyard/euler/}{Geometry Junkyard}}. To prove it, we can start by observing that any mesh of sphere topology has a continuous, bijective mapping onto the sphere. This entails that we can transform the geometry of such a mesh into that of a sphere while keeping the structure of the mesh unchanged. Furthermore, we can ensure that the edges of the original mesh map to \textit{great arcs} in the spherical mesh. \marginnote{\footnotesize\itshape A greeat arc is a segment of a circle formed as the intersection of the sphere with a plane that contains its centre.} With this in place, the polygons of \mesh map to spherical polygons which will be necessary. Figure~\ref{fig:armadillo-sphere} shows an example where a mesh of spherical topology is mapped onto a sphere. 

We begin by observing that the area of the (unit) sphere is $4\pi$. Thus, since the spherical polygons obtained from mapping onto the sphere must cover the sphere precisely, it is true that
\begin{equation}
    \sum_{f=1}^{\F} A_f = 4\pi \,,
    \label{eq:area-sphere}
\end{equation}
where $A_f$ is the spherical area of face $f$.

We know that the area of a spherical polygon is given by Harriot's formula,
\begin{equation}
    A_f = \sum_{i=1}^{k_f} \alpha_i - (k_f - 2) \pi \,,
    \label{eq:Harriot}
\end{equation}
where $k_f$ is the number of sides of $f$ and $\alpha_i$ is the angle of corner $i$ of the polygon. See Tristan Needham \cite{needham2021visual} for a proof of this formula for spherical triangles. The generalisation to spherical polygons is straightforward since polygons can always be triangulated. We can plug \eqref{eq:Harriot} into \eqref{eq:area-sphere}, obtaining
\begin{equation}
    \sum_{f=1}^{\F} \left(\sum_{i=1}^{k_f} \alpha_i - (k_f - 2) \pi \right)= 4\pi
\end{equation}

To proceed, we use the fact that the angles around each vertex sum to $2\pi$, though it may be worth pointing out that this is only true because the polygons have been transformed into spherical polygons. If the edges had been straight then the angle sum would be less than $2\pi$ since the corners would not have been flat in the limit. Now, we have
\begin{equation}
    2 \pi V - \sum_{f=1}^{\F} (k_f - 2) \pi = 4\pi 
\end{equation}

The final observation that we require is that the sum of the sides of the faces, $k_f$, is twice the number of edges, $E$, since every edge is shared by two faces. Reducing the formula, we get
\begin{equation}
    2 \pi V - 2\pi E + 2 \pi F = 4\pi     
\end{equation}
from which we can obtain Euler's formula just by dividing by $2 \pi$.
\section{Quantitative Topology}
Some topological properties are quantitative. In the context of polygonal meshes or surface more generally, we identify the numbers of 
\begin{itemize}
    \item shells (or connected components),
    \item holes (referred to as \textit{genus}), and
    \item boundary curves (which, confusingly, are also thought of as holes).
\end{itemize}
The number of shells is simply the number of separate pieces of the mesh, and the number of boundary curves is the number of distinct boundary cycles. Finally, the genus is the number of holes (in the sense of a hole in a donut) also sometimes referred to as the number of handles of the shape since any hole can be transformed into a handle without tearing or stitching the shape which means that holes and handles are the same from a topological point of view. Figure~\ref{fig:fertility} shows two objects both of genus four, but on the right it is much more obvious that the genus corresponds to the number of handles than on the left.

None of these three properties can be obtained directly from Euler's formula. Fortunately, there is a generalisation known as the Euler-Poincaré formula:
\begin{equation}
    V-E+F=2(s-g)-b = \chi \,,
    \label{eq:Euler-Poincare}
\end{equation}
where $s$ is the number of connected components (or shells), $g$ is the {genus}, $b$ is the number of boundary curves, and $\chi$ is known as the \textit{Euler characteristic}.

\paragraph{Shells}
$s$ is easy to interpret. It is simply the number of connected components of the mesh.
As, we recall there is no notion of solids when we deal with polygonal meshes. Hence, we do not distinguish between mesh components that represent interior voids and separate shapes even if it may be easy to tell these cases apart based on the geometry. If the mesh does have sphere topology, $s=1$, $g=0$, and $b=0$, and we immediately see that the formula reduces to the familiar one in this case. The generalisation to $s>1$ is trivial. If we have two objects of sphere topology and add up their vertices, edges, and faces, we get 
\begin{equation}
    V-E+F=2s \,.
\end{equation}

Proceeding to polygonal meshes with boundaries, we can assume that there is only a single component, $s=1$, without loss of generality.

\paragraph{Boundary Curves}
$b$ is the number of boundary curves. Let us first consider the simple case where a boundary is introduced by removing just a single face. It is clear that
\begin{equation}
V-E+F'=2 - 1
\end{equation}
where $F'=F-1$ since we removed one face. In other words, removing a single face just amounts to subtracting one from the right-hand side. What if the hole is bigger than a single face? We soon realise that removing more than one face does not change the left-hand side: when a face is removed to enlarge a hole, we also remove the edge that separates this faces from the hole. In some cases a face shares more than one edge with a hole, but then a vertex is also removed, and the end result would always be that the left-hand side is not changed further. We arrive at
\begin{equation}
V-E+F=2s-b \,.
\end{equation}


\paragraph{Genus}
Arguably, $g$, the genus of the shape, is the ``poster entity'' of topological invariants, and it is equal to the number of handles of the shape. A handle is the same as a hole (in the donut sense), since we can always deform (and homeomorphically map) a shape with $N$ holes into a sphere with $N$ handles. Thus, $g=1$ for a torus which has a single handle and $g=0$ for a sphere that has none.

\begin{figure}
    \includegraphics[width=0.49\textwidth]{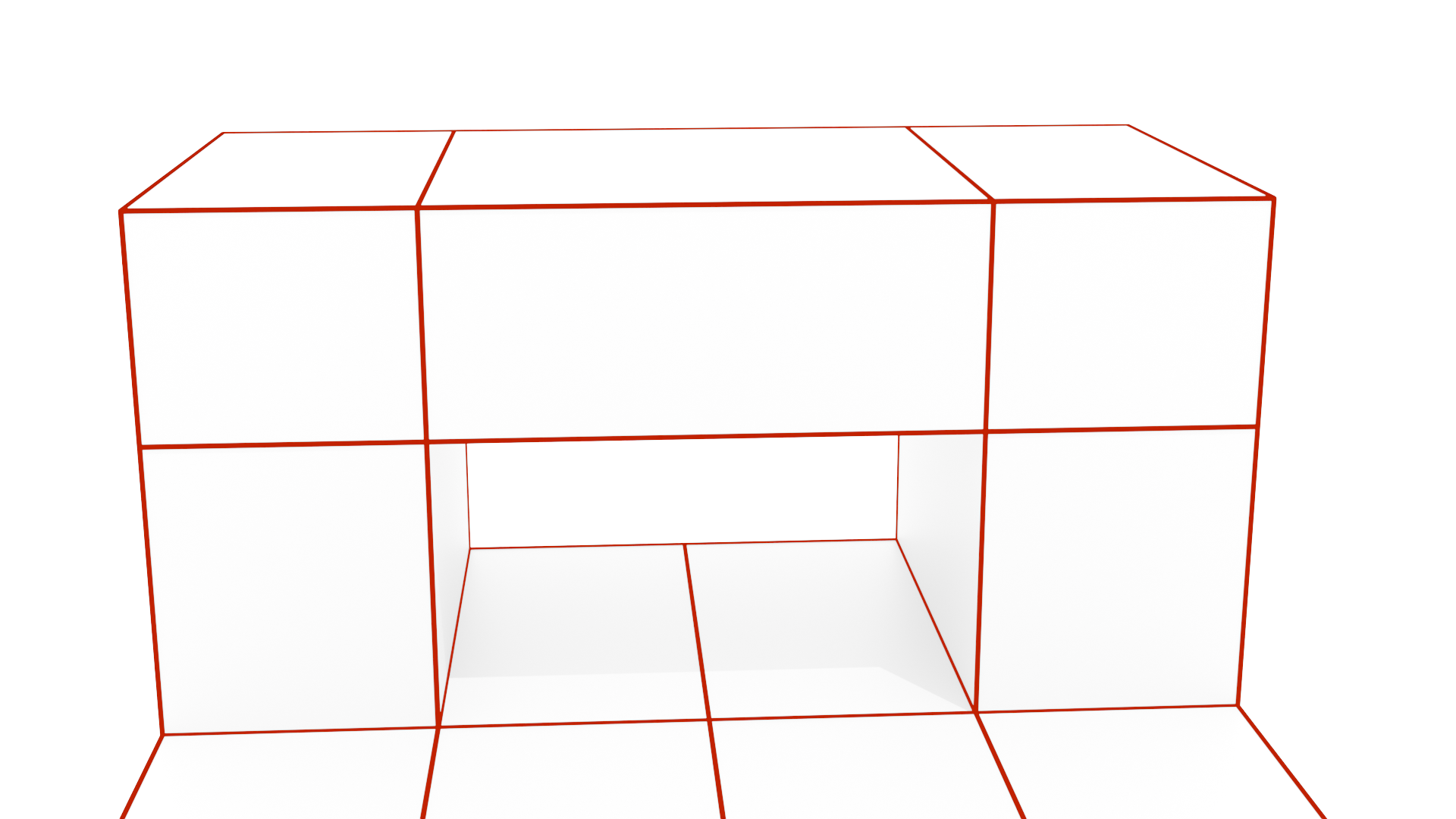}
    \includegraphics[width=0.49\textwidth]{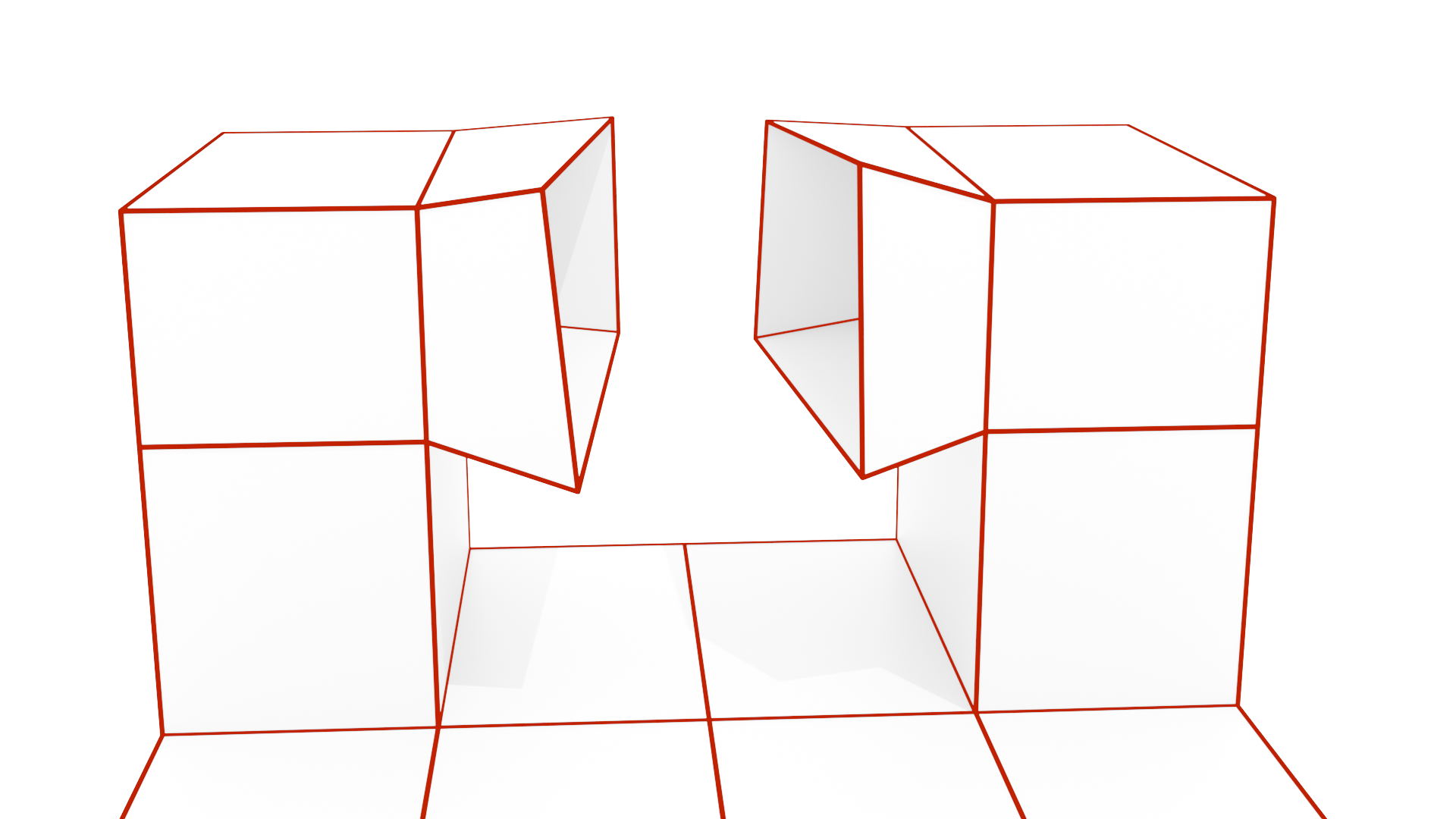}
    \caption{Left: A polygonal mesh with handle. Right: The same mesh after the handle was cut. It is easy to see how the geometry has proliferated due to the cut. There is an additional eight vertices, twelve edges and four faces. Moreover, the cut has introduced two holes in the mesh.}
    \label{fig:cut-handle}
\end{figure}

To understand how $g$ affects the Euler formula, we only have to consider what happens when we cut a handle \cite{courant1996mathematics}. The handle is part of the mesh, and cutting it we would  cut $N$ edges and slice through their incident $N$ faces. Where we cut each edge we need to insert two new vertices -- one on either side of the cut -- and the new vertices need to be connected by edges -- again an edge on either side. This means that we introduce $2N$ new vertices and $3N$ new edges since we split $N$ edges and need $2N$ edges to connect the new vertices. We also split $N$ faces introducing $N$ new faces. Putting this together we see that
\begin{equation}
V'-E'+F'= (V+2N) - (E + 3N) + (F+N) = V-E+F \,. 
\end{equation}
Yet something changed since we never closed the two holes where the handle was cut as shown in Figure~\ref{fig:cut-handle}. Thus, in terms of the Euler characteristic, a handle is equivalent to two holes. It must haven been true that $V-E+F=2-2=0$ \textit{before} we made the cut, if \mesh was a topological torus with a single handle. It follows that we have to reduce the right-hand side by two for every handle. This leads to 
\begin{equation}
V-E+F=2s - 2g - b \,.
\end{equation}
Thus, the general Euler-Poincaré formula has been proven -- or at least the proof has been sketched since we glossed over some details.

\subsection{Computing $s$, $g$, and $b$}
\label{sec:computing-sgb}
Fortunately, it is easy to find the values of $s$, $g$, and $b$ if we are given a polygonal mesh, \mesh. Finding $V$, $E$, and $F$ on the left hand side of \eqref{eq:Euler-Poincare} is just a matter of counting. However, we must remember to only count edges and vertices once when they are shared by several faces of \mesh. If the faces of \mesh are stored in such a way that we cannot tell whether edges or vertices are shared, it is necessary to first match up corners and polygon edges that correspond to identical vertices and edges of \mesh. 

$s$ can be found using ``flood filling''. Start from any given vertex and visit adjacent, unvisited vertices recursively. When this procedure terminates, pick an unvisited vertex and repeat. $s$ is the number of times this process is started. 

$b$ can be found in a very similar fashion. We can identify boundary edges because they only have a single incident face. Starting from any boundary edge, we visit all the connected boundary edges, tracing out its cycle (c.f. Section~\ref{sec:cycles}). Once a boundary cycle has been found, we mark all its edges as visited and pick the next unvisited edge. $b$ is the number of times this process was started.

$g$ is now the only term that is missing and we can easily solve for $g$ using \eqref{eq:Euler-Poincare}.
\section{The Betti Numbers}
Having proven that the Euler-Poincaré formula is true does not seem to have been very helpful for understanding it. However, there is a clue in its structure: it looks like the beginning of an alternating sum since the first term, $V$, is added, the next term, $E$, is subtracted, and the final term, $F$, is again added. We will now show that $V-E+F$ is equal to another alternating sum which has a much more clear interpretation.
\subsection{Cycle Instigators}
\label{sec:instigators}
The first task is to count how many vertices, edges, and faces, we can choose without creating any cycles. To this end, we will introduce the notion of a cycle \textit{instigator}. Given a set of mesh edges, a cycle instigator is an edge, $e^i$, that becomes part of a cycle when added to the set. To be precise, $e^i$ is a cycle instigator with respect to a set of edges $\mathcal{T}$ iff
\begin{equation}
\exists \mathcal C: \mathcal C \subseteq \mathcal{T} \cup e^i \wedge e^i \in \mathcal C
\end{equation}
where $\mathcal C$ is a cycle. Now, the question we ask is how many edges from the total set of edges, $\mathcal E$, we can add to $\mathcal T$ without breaking the property that $\mathcal T$ contains no cycles? Put differently, what is the maximal number of edges we can add without including an instigator?

\begin{figure}
    \centering
    \includegraphics[width=0.9\textwidth]{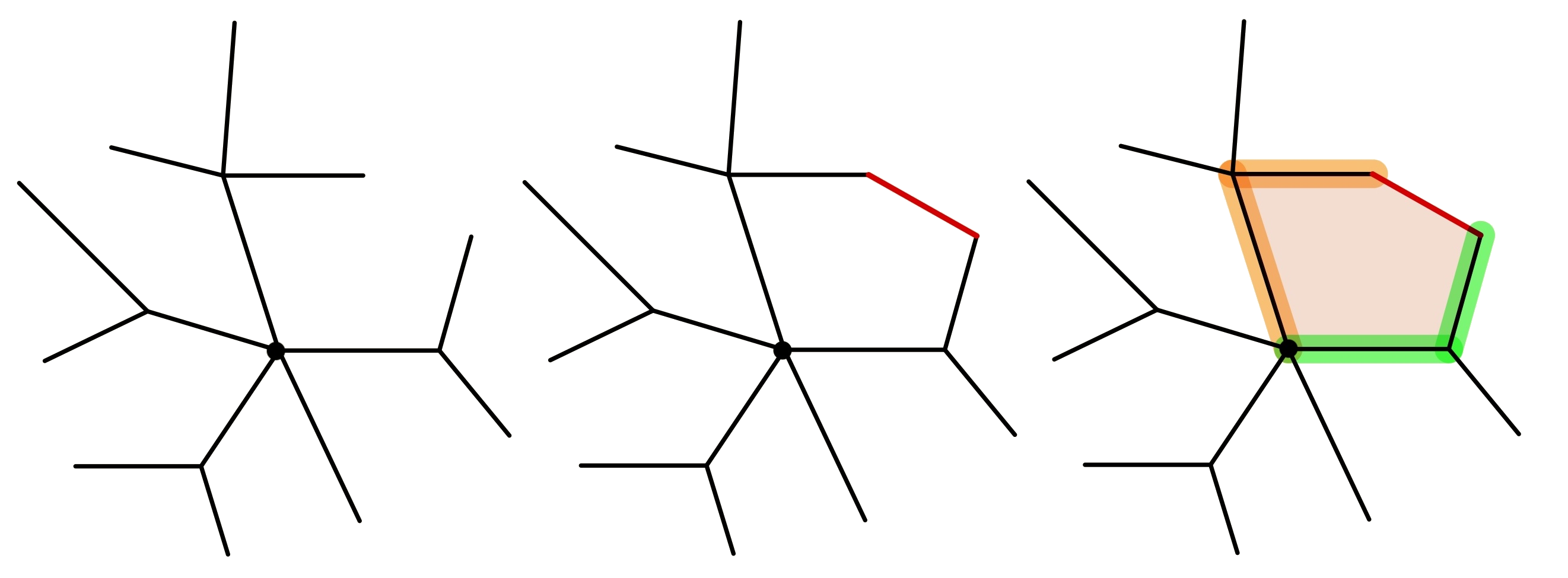}
    \caption{A tree is shown on the left with a black dot that indicates the root. In the middle image, a cycle-instigator has been added (red), and on the right the orange and green curves connect the end points of the instigator edge with the root. Adjoining the instigator edge with the orange and green edges, we clearly form a cycle (that bounds the pink region).}
    \label{fig:cycle-tracing}
\end{figure}
In fact, this question is easy to answer. Assuming \mesh is a single component, a \textit{spanning tree} of the vertices of \mesh is a collection of edges which is incident on every vertex of $\V$. Now, a tree contains no cycles by definition, and we can create a tree by first adding any edge of the mesh to $\mathcal T$ and subsequently adding only edges that connect a vertex already incident on $\mathcal T$ to a vertex on which $\mathcal T$ is not incident. If we count the edges while adding them, it is clear that our count reaches $V-1$ since we add one edge for every vertex except the first. To answer the question, we now note that any edge that does not belong to the tree defines a unique cycle and is thus an instigator. We can obtain the instigated cycle by adjoining the instigating edge to the path we obtain by following the tree edges from its end points back towards the root until they meet at a vertex (which may be the root vertex). This procedure is shown in Figure~\ref{fig:cycle-tracing}. We conclude that the maximal number of non-instigating edges is the number of edges in a spanning tree, i.e. $V-1$, and we can now define the number of non-instigating edges, $E_N$, and the number of cycle-instigating edges, $E_C$, as follows:
\begin{eqnarray}
\label{eq:e_n}
E_N &=& V - 1 \, ,\\
E_C &=& E - E_N = E - V + 1 \,.
\end{eqnarray}
Note that in the case of multiple connected components 
\begin{eqnarray}
    \label{eq:e_n_multiple}
    E_N &=& V-s\,,\\
    E_C &=& E - E_N = E - V + s \,.
\end{eqnarray} 
since we have to construct a tree for each component, $i$, and the tree for each component consists of $V_i - 1$ edges where $\sum_{i\in[1..s]} V_i = V$.

In the case of vertices, the task of dividing up the vertices is trivial. All vertices are cycles, hence all vertices instigate cycles, and
\begin{eqnarray}
V_N &=& 0 \, ,\\
V_C &=& V \,.
\end{eqnarray}
For faces, the situation is also straightforward. A mesh is a cycle only if it is closed. This means that we can pick all but one of the faces without creating a cycle. In other words, if the mesh is watertight, the maximum (and minimum) number of non-instigating faces is the number of faces minus 1:
\begin{eqnarray}
F_N &=& F - 1\, ,  \\
F_C &=& 1
\end{eqnarray}
If the mesh is not watertight, $F_N=F$ and $F_C=0$. When we are dealing with meshes that consist of multiple components which may or may not be watertight, then
\begin{eqnarray}
F_N &=& F - s_w\, ,  \\
F_C &=& s_w \, ,
\end{eqnarray}
where $s_w$ is the number of watertight components.

In all cases, the number of instigators and non-instigators sum to the original number of entities: $V=V_N+V_C$, $E=E_N+E_C$, and $F=F_N+F_C$. 
\subsection{Reformulating Euler's Formula}
With the division of mesh entities into instigators and non-instigators of cycles, we can express the right-hand side of the Euler-Poincaré formula slightly differently,
\begin{equation}
    V-E+F = (V_N+V_C)-(E_N+E_C)+(F_N+F_C) \,,
\end{equation}
and rearrange the terms by shifting the non-instigating terms to the left,
\begin{equation}
    V-E+F = (V_C-E_N)-(E_C-F_N)+(F_C-0) \,,
\end{equation}
since $V_N=0$ and there is no non-instigating term of higher dimension than that for faces. Let us rename the three terms on the right-hand side:
\begin{eqnarray}
    \label{eq:betti-0}
\beta_0 &=& V_C-E_N\,,\\
    \label{eq:betti-1}
\beta_1 &=& E_C-F_N\,,\\
    \label{eq:betti-2}
\beta_2 &=& F_C-0 \,.
\end{eqnarray}
We now obtain
\begin{equation}
    \label{eq:euler-betti}
    V-E+F= \beta_0-\beta_1+\beta_2 \, .
\end{equation}
The left-hand side of \eqref{eq:euler-betti} is an alternating sum of mesh elements that increase in dimension. The right-hand side is a different alternating sum where each term is one of the \textit{Betti numbers}. Each Betti number is thus the difference between the number of cycle instigators and the non instigators one dimension higher.
\paragraph{Understanding the Betti Numbers}
The alternating sum $\beta_0 - \beta_1 + \beta_2$ explains more about the topology of the mesh than $V-E+F$ although the precise significance of the $\beta_i$ might still be a bit hazy. Below, we will try to clarify the notion of Betti numbers.

Consider a mesh with just a single vertex. That vertex is the sole entity of the mesh and a cycle instigator. It is clear that $\beta_0 = 1 -0= 1$ and $\beta_1=\beta_2=0$. Add another vertex, and we can create cycles consisting of either vertex or both combined: the vertices are basis vectors in a space of {vertex cycles} and $\beta_0=2-0=2$ means that the dimension of the space of vertex cycles is two. Add a non-instigating edge that connects the two vertices, and we are back to the original situation, $\beta_0 = 2 -1 = 1$ since the now connected vertices no longer represent distinct vertex cycles. In other words, adding the (non-instigating) edge collapses the dimension of the space back down to one.

Likewise, when we add an instigating edge, it is not a single cycle but a subspace in a space of edge cycles that is created. This is because the instigator defines a unique cycle that can be construed as a basis vector. If we later add a non-instigating face that fills in the cycle, the boundary of that face also represents a basis vector which we use to project out a dimension of the space, decreasing $\beta_1$ by 1. A geometric way to think about it is that adding the face allows the curve to ``slip across it'' making edge cycles indistinguishable regardless of which way they go around the face. The upshot is that $\beta_1$ represents the dimension of the space of edge cycles that do not bound faces. 

The analysis of $\beta_2$ is simple. There is only one cycle-instigating face since a polygonal mesh is a face cycle only when it is watertight. When this face is added, $\beta_2=1$. If polygonal meshes were considered solids when watertight, the next step would be to add a non-instigating solid that fills in the face cycle, reducing $\beta_2$ to zero. Thus, the fact that we consider polygonal meshes to be hollow explains why the Euler characteristic, $\chi$, of a topological sphere is 2: a topological sphere is a single connected component with a single shell and we add one for each of these two features.      
\paragraph{Incremental Construction}
In the more general setting of simplicial complexes, Delfinado and Edelsbrunner \cite{delfinado1995incremental} describe an incremental algorithm for computing Betti numbers. Instead of splitting mesh elements into instigators and non-instigators up front, they rely on a so-called \textit{filtration}. A filtration is an arbitrary ordering of all elements of the mesh with the condition imposed that a tetrahedron can only appear in the ordering after all its triangular faces, a face can only appear after all of its edges and an edge only after its vertices. 

For each element of a mesh, \mesh, we will use $k$ to represent its dimension: $k=0$ for vertices, $k=1$ for edges, $k=2$ for faces, and $k=3$ for tetrahedra. We will let $\mathcal M_i$ denote the mesh where all elements up to and including the $i$-th element have been added.

\mesh is now built by adding the elements of the filtration one at a time. For each element, the $k$-th Betti number where $k$ is the dimension of the element, is increased by one if it is an instigator (in our terminology) with respect to $\mathcal M_i$. If it is a non-instigator, the $(k-1)$-th Betti number is decreased by one. Although the elements are added and subtracted one at a time, this leads to the same result as Equations \eqref{eq:betti-0}, \eqref{eq:betti-1}, and \eqref{eq:betti-2}. However, it is necessary to check after each insertion whether a cycle has been created.
\section{Examples}
\begin{figure}
    \centering
    \includegraphics[width=\textwidth]{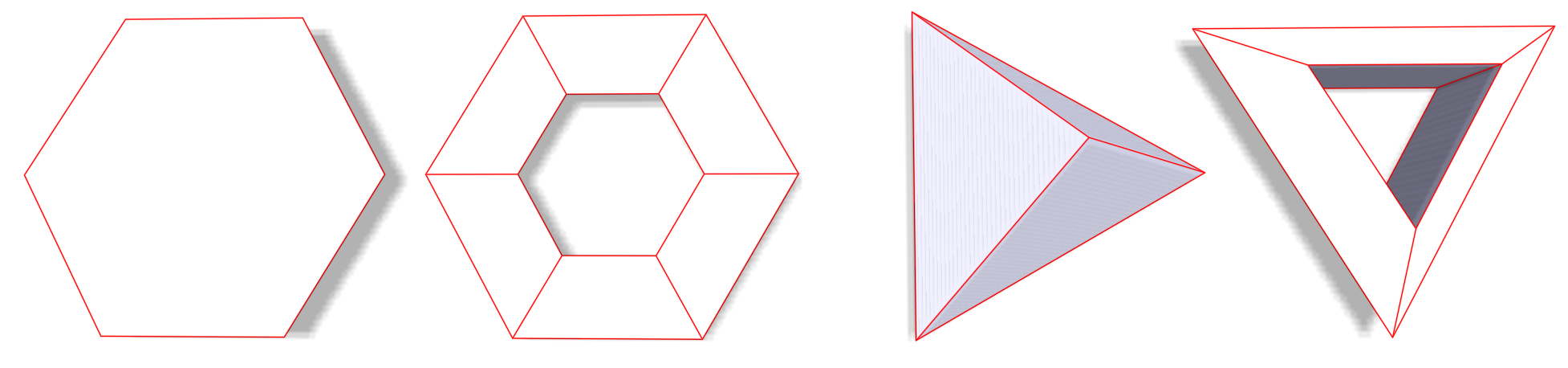}
    \caption{From left to right this image shows a hexagon (topological disc), annulus, tetrahedron (topological sphere), and torus. While all of these objects are simple they have different topology.}
    \label{fig:simple-shapes}
\end{figure}

With the formulas that we have derived, it is straightforward to compute the quantitative topological properties of polygonal meshes. Figure~\ref{fig:simple-shapes} shows four shapes which are topologically quite different.

To find $s$, $g$, and $b$ we follow the procedure outlined in Section~\ref{sec:computing-sgb} and the Betti numbers are simply computed using Equations \eqref{eq:betti-0}, \eqref{eq:betti-1}, and \eqref{eq:betti-2}. The quantitative properties of the four example shapes are summarised in Table~\ref{tab:quantitative-properties}.

\marginnote{\footnotesize\itshape \includegraphics[width=4.2cm]{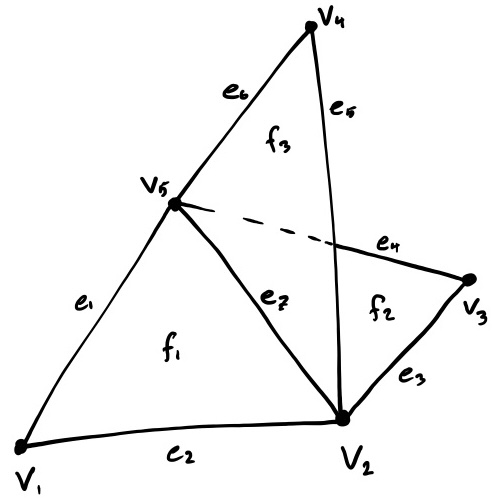}}
The two triples of numbers usually do not agree, except that $\beta_0 = s$, but since $V-E+F=2(s-g)-b=\beta_0-\beta_1+\beta_2$ both of these quantitative descriptions of the topology of an object do contain the same information; indeed it seems as if the information is already contained in the numbers $V$, $E$, and $F$.

Of course, we must be a little careful. Imagine an object which is a single edge shared by three triangles (shown right). This object consists of 5 vertices, 7 edges, and 3 faces. We cannot compute the number of boundary curves, since they do not form closed cycles. We can compute the Betti numbers, though, and get $\beta_0=1$, $\beta_1=0$, and $\beta_2=0$. In other words, while this object is non-manifold it has the Betti numbers of a topological disc. Likewise, we could turn the annulus into a Möbius band and the torus into a Klein bottle -- all without changing their respective Betti numbers since these objects have the same numbers of vertices, edges, and faces.

\marginnote{\footnotesize\itshape Note: some types of simplicial homology do allow us to distinguish between orientable and non-orientable surfaces \cite{giblin2013graphs}.}
This informs us that the quantitative topological properties should be understood in light of the qualitative topological properties. It is not enough to know $V$, $E$, and $F$. We also need to know if the polygonal mesh is manifold, whether it has a boundary, and if it is orientable -- in order to meaningfully compute the quantitative topological properties.
\begin{table}
    \centering
    \caption{Example meshes for which we have computed the numbers of non-instigating and instigating elements, $s$, $g$, $b$, and $\beta_i$.}
    \label{tab:quantitative-properties}
    \begin{tabular}{r|llllll|lll|lll}
        Mesh & $V_N$ & $V_C$ & $E_N$ & $E_C$ & $F_N$ & $F_C$ & $s$ & $g$ & $b$ & $\beta_0$ & $\beta_1$ & $\beta_2$ \\
        \hline
        Hexagon     & 0 &  6 &  5 &  1 & 1 & 0 &   1 & 0 & 1 & 1 & 0 & 0\\
        Annulus     & 0 & 12 & 11 &  7 & 6 & 0 &   1 & 0 & 2 & 1 & 1 & 0\\
        Tetrahedron & 0 &  4 &  3 &  3 & 3 & 1 &   1 & 0 & 0 & 1 & 0 & 1\\
        Torus       & 0 &  9 &  8 & 10 & 8 & 1 &   1 & 1 & 0 & 1 & 2 & 1\\
    \end{tabular}
\end{table}
\section{Computing a Cut Graph}
For some applications it is necessary to transform a mesh into a topological disc. If the mesh happens to be a topological sphere a single open cut suffices since we only need to ``puncture'' the sphere in order to turn it into a topological disc.

For shapes of higher genus, i.e. $g>0$, we need $2g$ closed cuts. Recall that for an object which consists of a single component and which has no boundary curves $\beta_1 = 2g$ and, for simplicity, we will assume that the object has no boundary curves. We will also assume there is just a single connected component. This means that the number of cuts we need to make is exactly the number of unique edge cycles which do not bound sets of faces.

In the following, we will present an algorithm originally due to Eppstein \cite{eppstein2003dynamic} and improved by Erickson and Whittlesey \cite{erickson2005greedy}. Later Keenan Crane presented the algorithm with some modifications \cite{crane_loops}.

\begin{figure}
    \centering
    \includegraphics[width=0.42\textwidth]{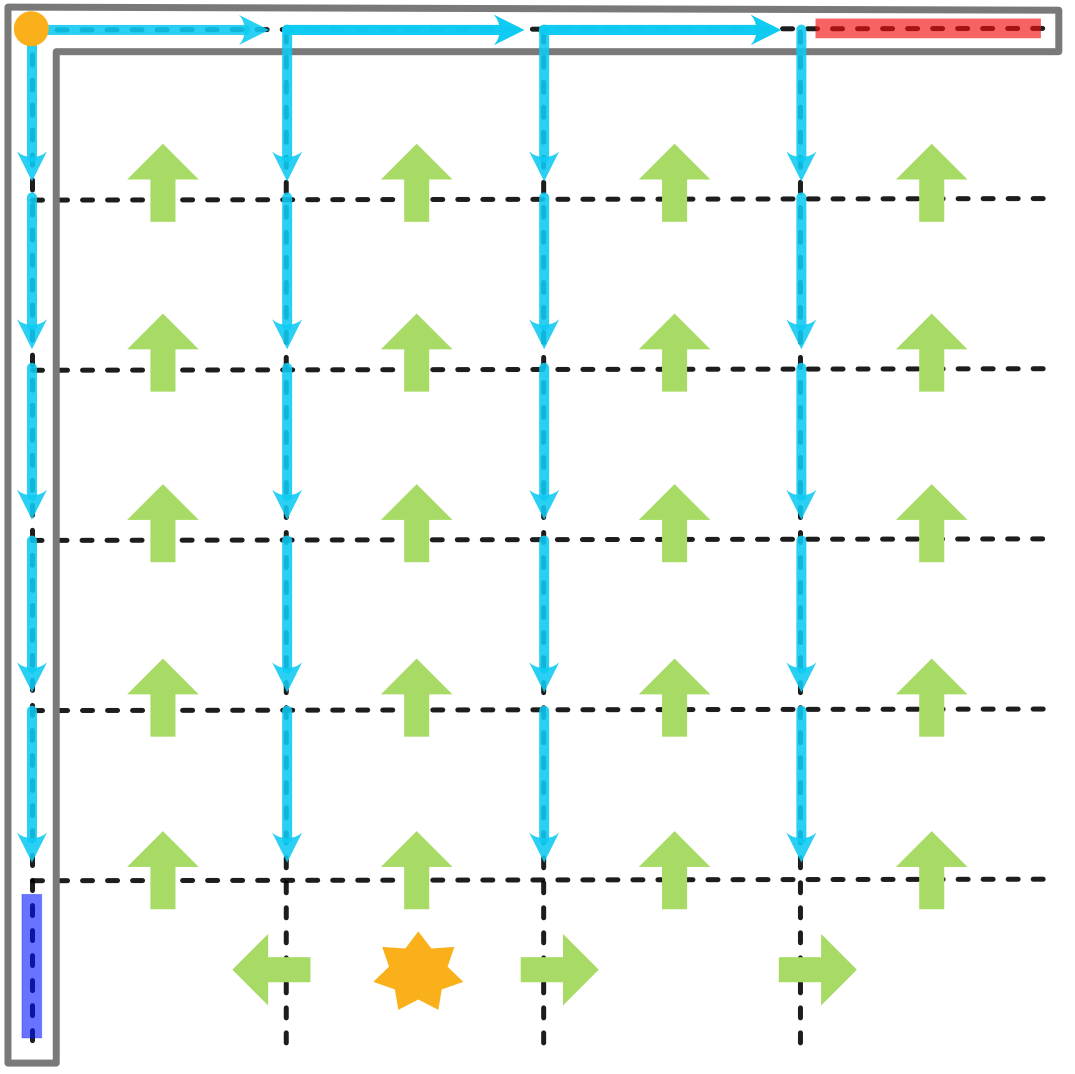}
    \includegraphics[width=0.57\textwidth]{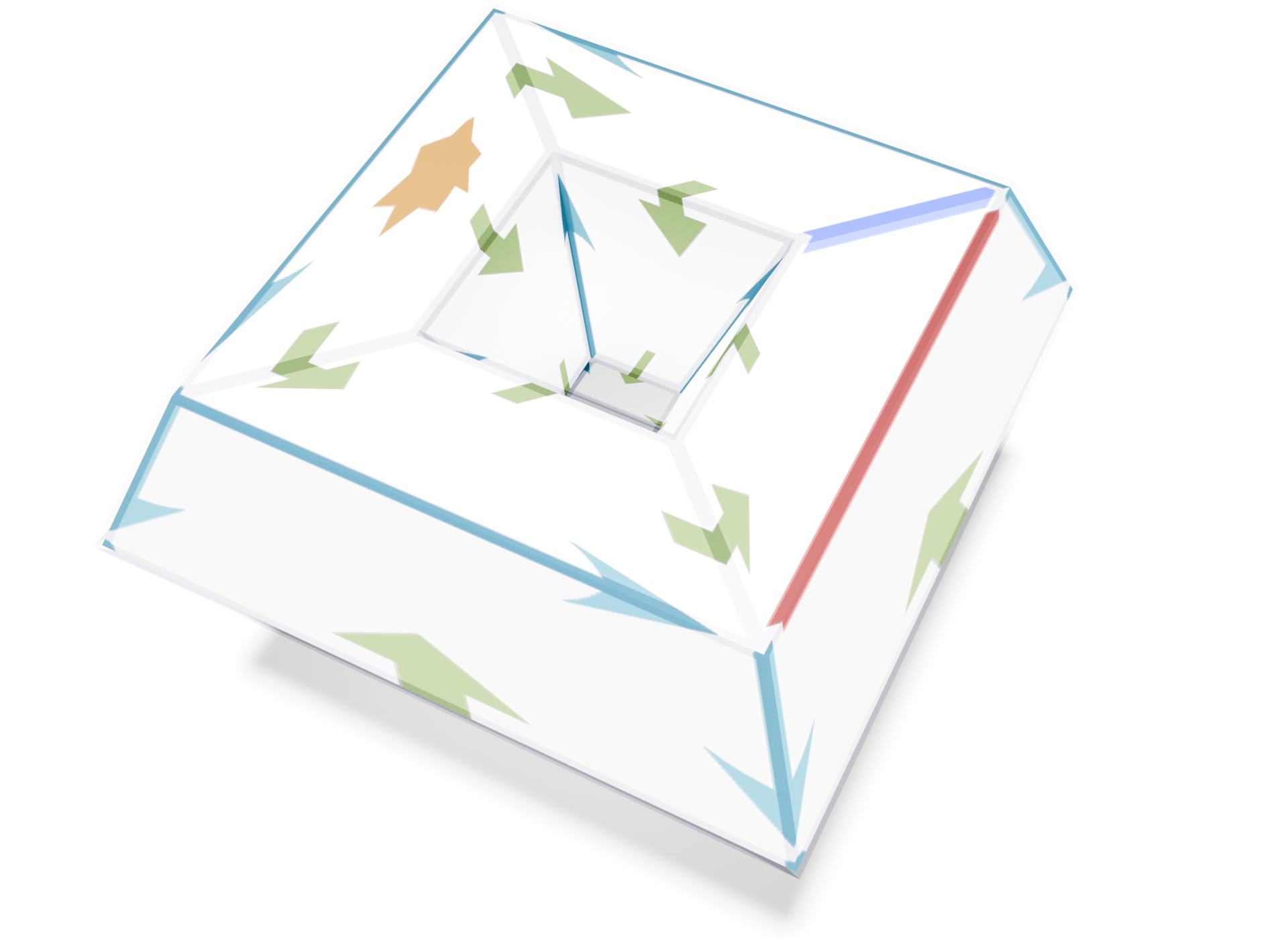}
    \caption{Left: A torus which has been cut and flattened to a disc. Note that we should consider the bottom to be attached to the top and likewise the right side is attached to the left. The edges of the shortest-path tree ($\mathcal T$) are shown pale blue, and the co-tree edges ($\mathcal T^c$) are those which are crossed by the bold green arrows. Edges point from parent to child in both trees, and the root of the shortest-path tree is the orange disc, and the root of the co-tree is the yellow star. The two edges highlighted by red and purple, respectively, are the cycle-instigating edges. The $\Gamma$ shaped frame highlights the edges in $\mathcal B$. Right: a rendering of the torus as a 3D object with the image on the left shown as a texture.}
    \label{fig:torus}
\end{figure}

\paragraph{Shortest-Path Tree Construction}
\marginnote{\footnotesize\itshape See Section~\ref{sec:instigators} for a definition of a spanning tree. Note also that a shortest-path tree should not be confused with a minimum spanning tree where the sum of the edge lengths is minimal. The shortest-path tree turns out to be advantageous \cite{erickson2005greedy}.}
Given a mesh, \mesh, the first step is to create a spanning tree, $\mathcal T \subset \E$, of the vertices. We can do this using Dijkstra's algorithm \cite{cormen2022introduction} which produces a tree of shortest paths from the root vertex to any vertex of the mesh. Just like any other spanning tree, the shortest-path tree contains $E_N = V-1$ edges. A shortest-path tree computed for a torus is shown in Figure~\ref{fig:torus}.

\paragraph{Co-Tree Construction}
The next step is to construct a so-called \textit{co-tree}, $\mathcal T^c \subset \E$, which is built from the dual of the mesh: a graph whose vertices are the faces of the mesh and two dual vertices are connected if the faces share an edge. We can build a co-tree using an algorithm that is similar to Dijkstra's. The key component is a queue that keeps track of faces which have not yet been visited but are neighbours of visited faces. 

An arbitrary starting face, $f_0$, is added to a queue. In each iteration, we pick a face $f$ from the queue and consider each neighbour, $f_n$, which shares an edge, $e$, with $f$. If $f_n$ has not been visited before and is not the parent of $f$ in the co-tree \marginnote{\footnotesize\itshape if we add $e$ to $\mathcal T^c$, $f$ will be the parent of $f_n$. Hence, $f_n$ must not be the parent of $f$.} and $e$ is not in the shortest-path tree, i.e. $e \notin \mathcal T$, 
then we add $f_n$ to the queue and mark it as visited while $e$ is added to $\mathcal T^c$. This process is iterated until the queue is empty. A co-tree for a torus is shown in Figure~\ref{fig:torus}.

Observe that the set of traversed edges must contain $F-1$ edges since the co-tree is a spanning tree over the faces. If we subtract the edges of the shortest-path tree from the total set of edges, we are down to $E_C=E-E_N$ edges. Now, we subtract $F_N=F-1$ edges, we have precisely $E_C-F_N=\beta_1$ edges left.
\paragraph{Constructing the Cut Graph}
The cut graph is used to turn the mesh into one that has disc topology. However, there is actually a simple procedure that produces a mesh of disc topology, $\mathcal M'$, without any need for the cut graph. This procedure is to cut along every edge except those that belong to $\mathcal T^c$. Since the co-tree is connected and has no loops, $\mathcal M'$ has disc topology.

The downside of this procedure is that far more edges than necessary are cut. Instead, we might opt for the following procedure. Trace back from the end points of each of the $\beta_1$ cycle-instigating edges along the edges of the shortest-path tree \marginnote{\footnotesize\itshape Tracing a loop was illustrated in Figure~\ref{fig:cycle-tracing}} thereby producing $\beta_1$ cycles. Let the set of edges in this set of cycles be denoted $\mathcal B \subset \E$. The \textit{cut graph} is now defined as the graph consisting of the edges in $\mathcal B$ (shown framed in Figure~\ref{fig:torus}).

We now claim that cutting \mesh along the cut graph, $\mathcal B$, turns it into a topological disc which we will denote $\mathcal M^*$.

To prove this, observe that if we were to cut $\mathcal M^*$ along all not-already-cut edges that belong to $\mathcal T$, i.e. edges in $\mathcal T \backslash \mathcal B$, we arrive at  $\mathcal M'$ which is known to have disc topology. Since $\mathcal T \backslash \mathcal B \subset \mathcal T$, it does not contain any cycles. In other words, we can transform $\mathcal M^*$ into $\mathcal M'$ without additional cyclic cuts. Hence $\mathcal M^*$ must be a topological disc -- possibly with additional boundary curves. 

We state without proof that if the cut graph consists of a single connected component, a disc topology mesh is produced; otherwise a disc with $n-1$ holes where $n$ is the number of connected components of the cut graph.

\subsection{Connections to Algebraic Topology}
The edge-cycles that belong to the cut graph are known as \textit{generators} of the \textit{classes} in the first homology group. By adding boundaries of collections of faces to generators and combinations of generators, it is possible to obtain any edge cycle on a mesh. The generators produced by the described algorithm are also sometimes called a \textit{homotopy basis}, e.g. by Keenan Crane \cite{crane_loops}. Grossly simplifying, the difference between a homotopy basis and a homology basis is that the generators of a homotopy basis have a common base point (the root of our shortest-path tree). For a homology basis we could choose different generator curves which do not have the same starting point. This can be advantageous if we care about the geometry of the generators. For instance, the HanTun Algorithm by Dey et al. \cite{dey2008computing} finds handle and tunnel loops which together form a homology basis. Note, though, that in this case the generators of the homology basis may not form a connected cut graph.
\subsection{Examples}
With the algorithm just described, we can compute a cut graph that turns a closed mesh of any genus into a topological disc.

For a torus, $\beta_1=2$ and the two instigator edges are shown in purple and red in Figure~\ref{fig:torus}. Note that in case of the red edge, the right end point is already connected to the root of the shortest path tree; in order to close the loop, we just have to trace back along the horizontal edges to the root vertex (in the top left corner). Likewise, to close the loop defined by the purple instigator, we simply trace vertically from the top vertex of the instigator to the root. It is clear that cutting along these two curves will turn the torus into a rectangle. Note that the top and bottom edges are really the same edges, and the same is true of the right side and left side edges.

\begin{figure}[p]
    \centering
    \includegraphics[width=0.29\textwidth]{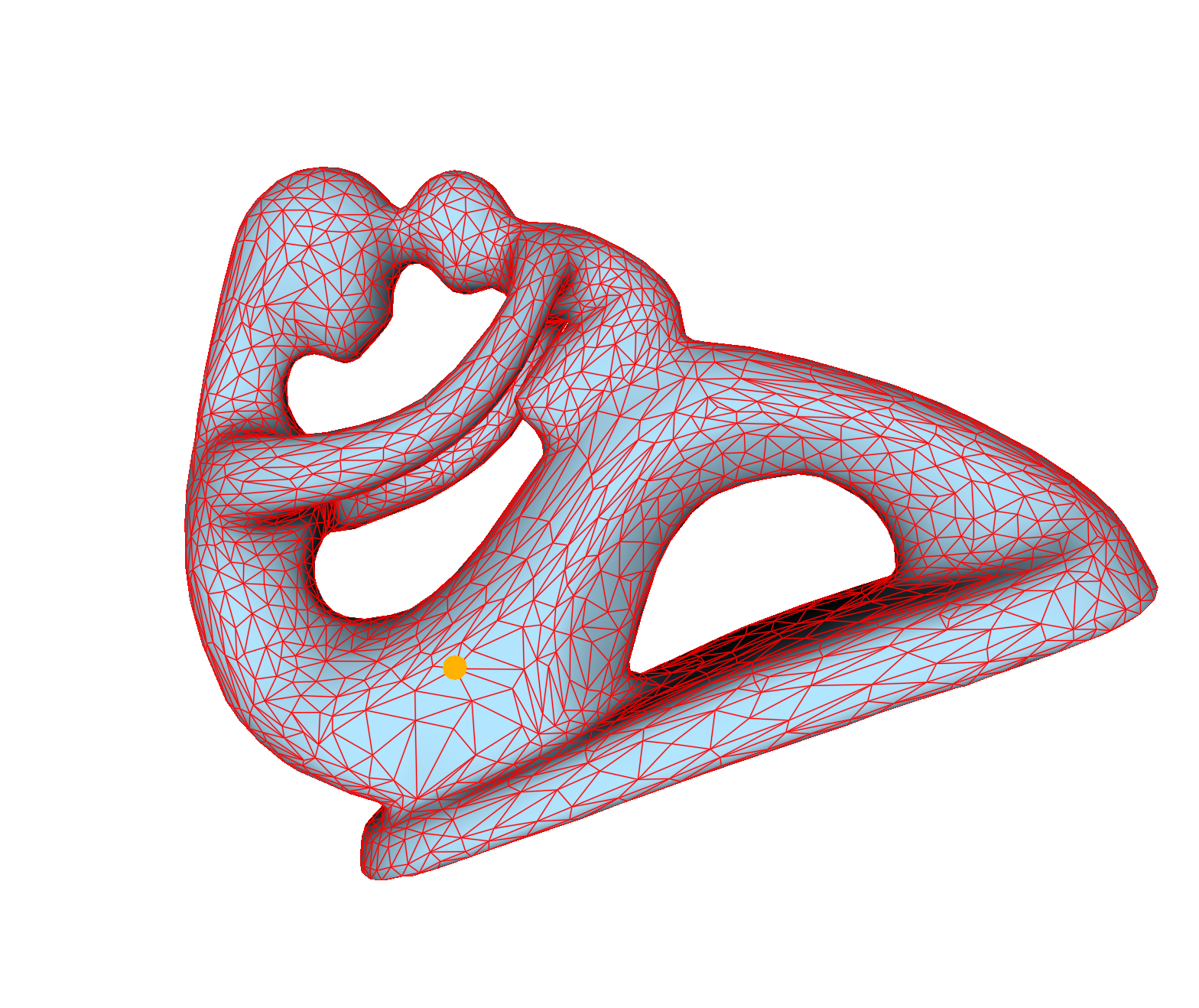}
    \includegraphics[width=0.29\textwidth]{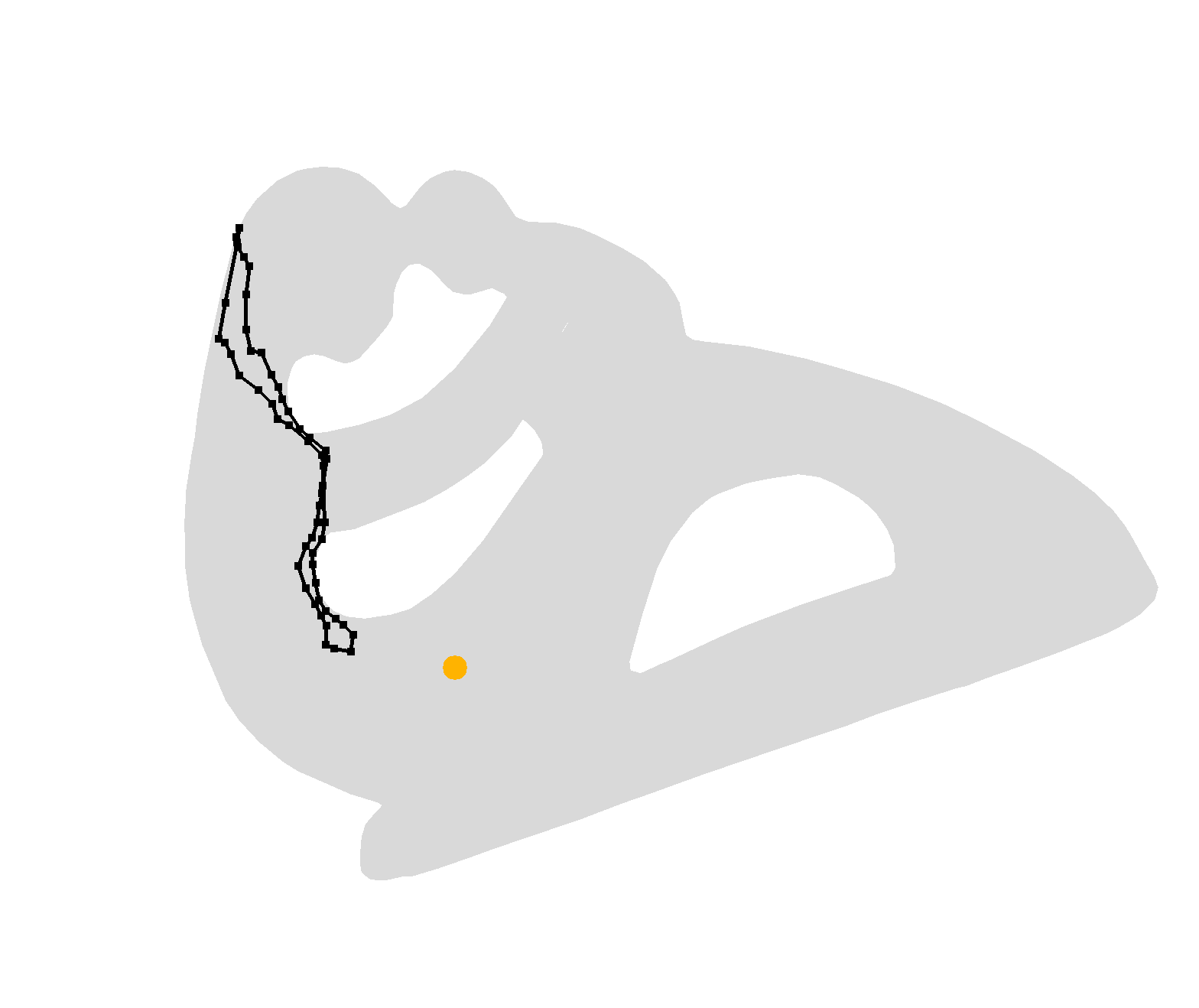}
    \includegraphics[width=0.29\textwidth]{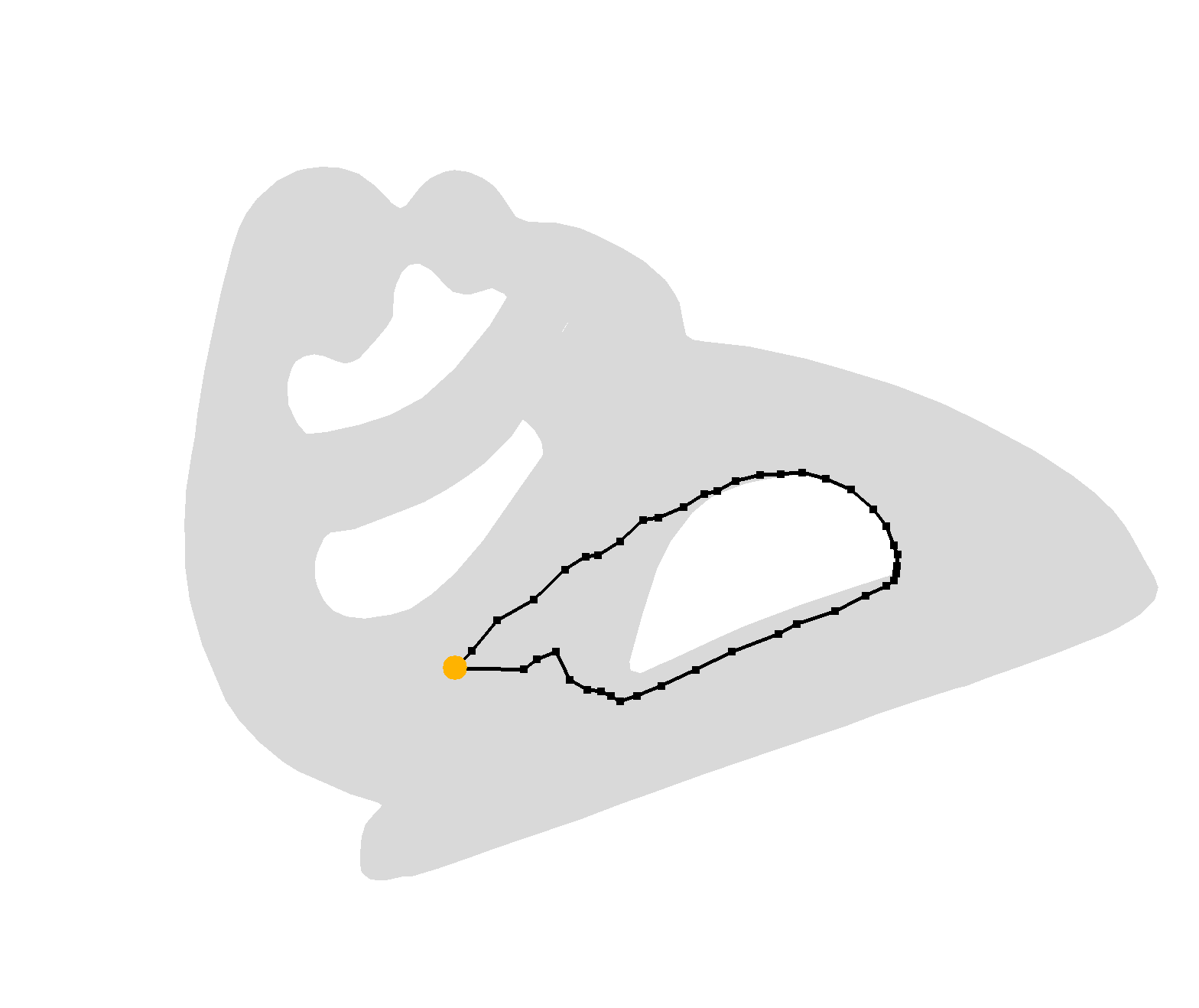}\\
    \includegraphics[width=0.29\textwidth]{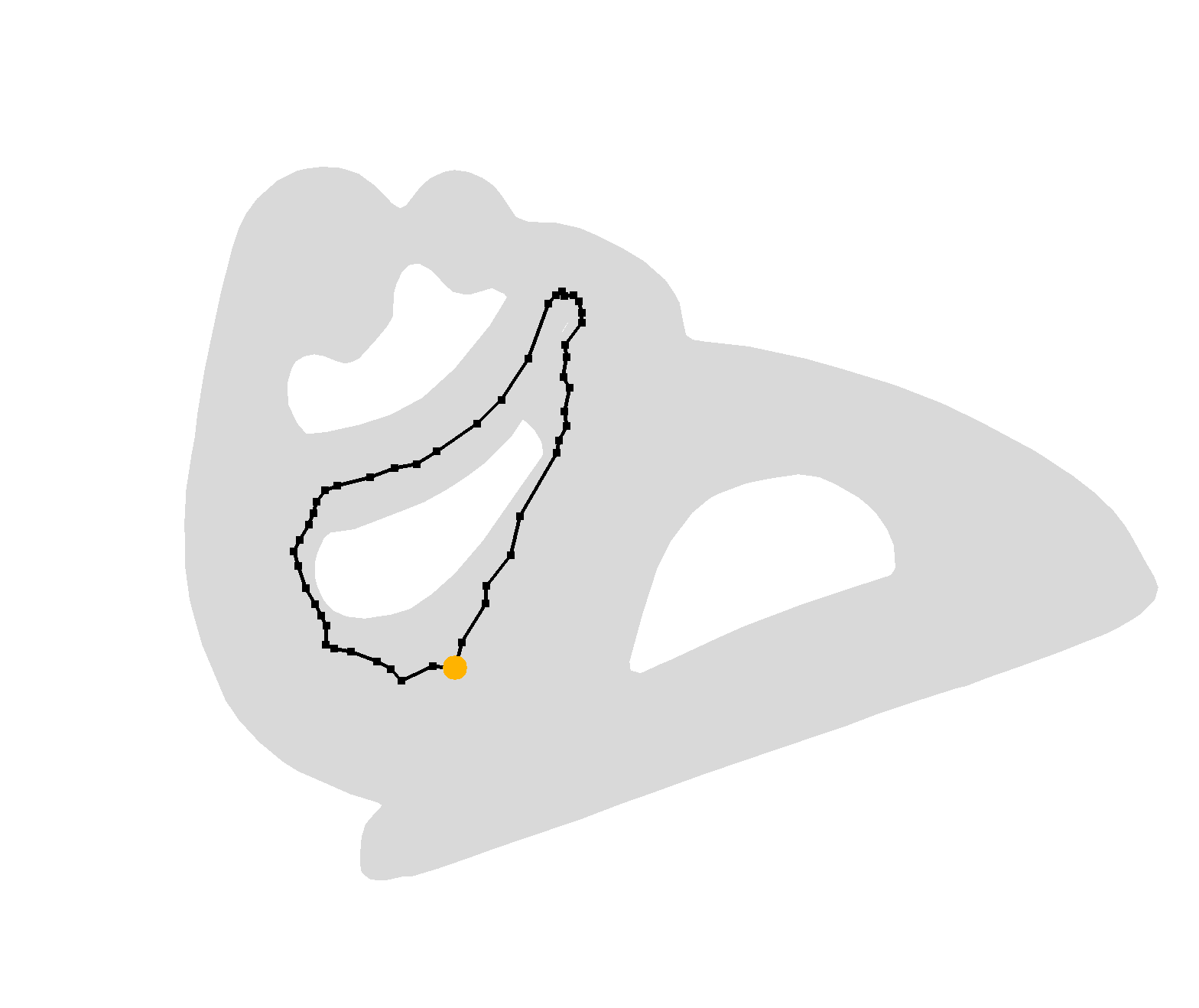}
    \includegraphics[width=0.29\textwidth]{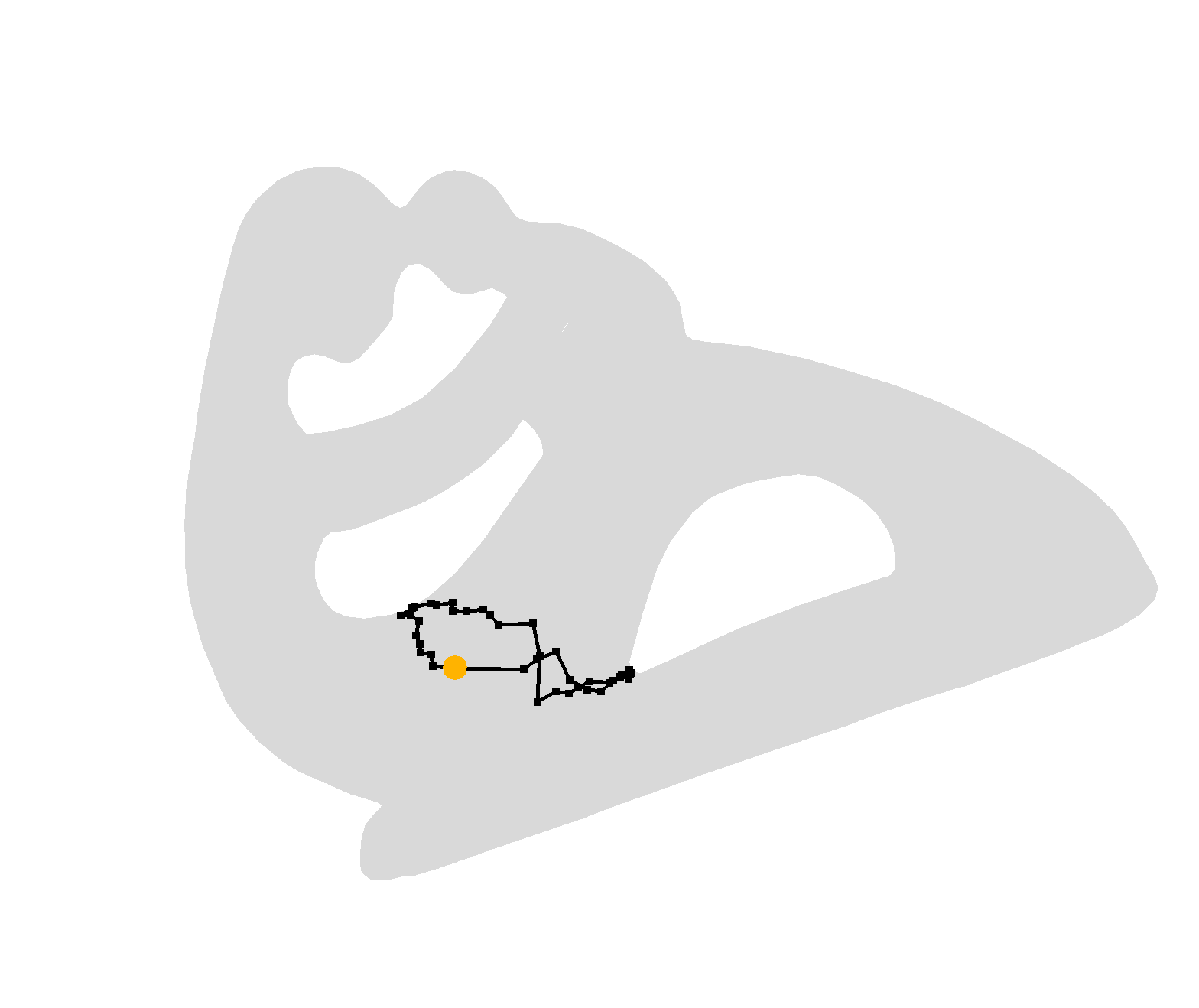}
    \includegraphics[width=0.29\textwidth]{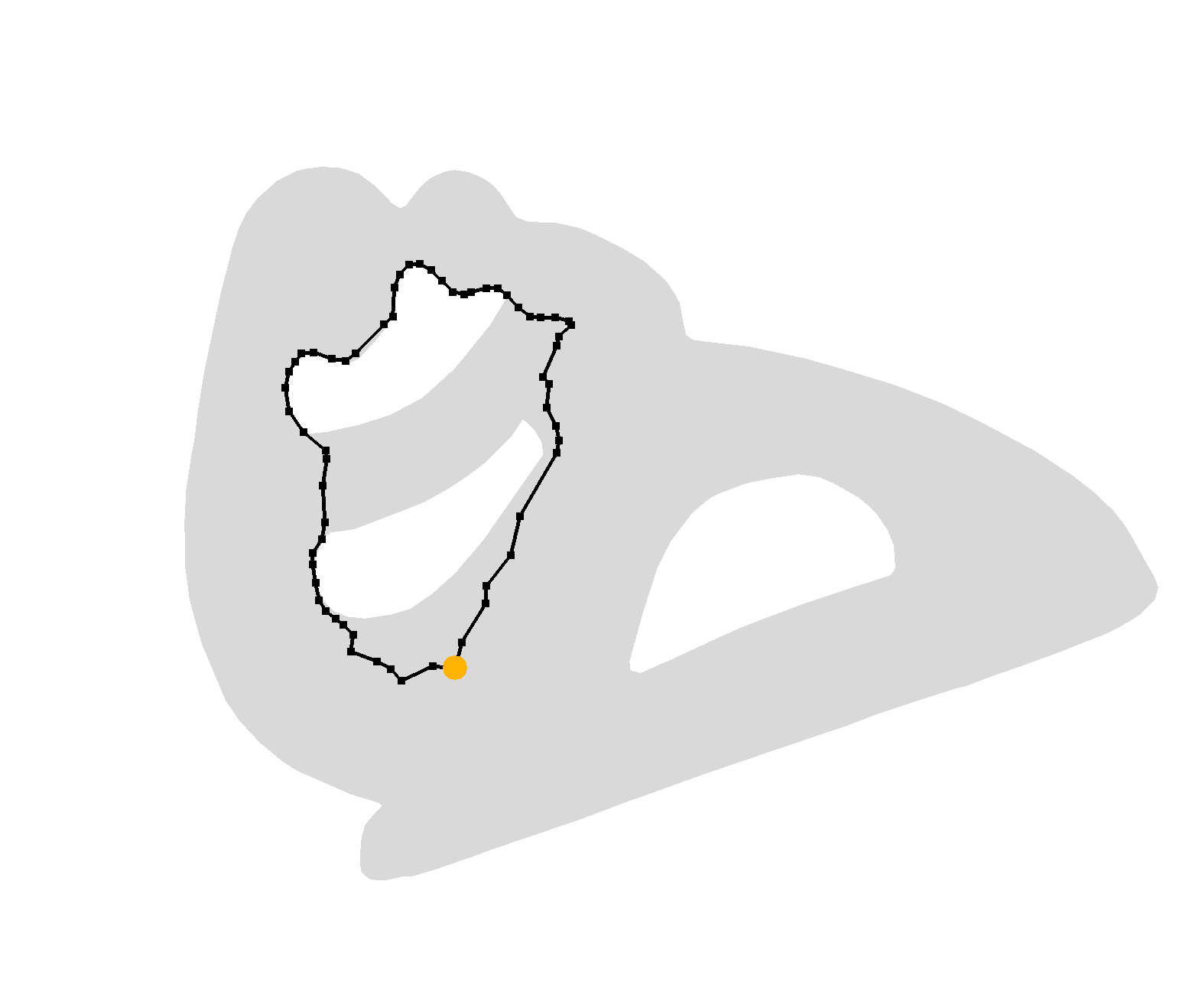}\\
    \includegraphics[width=0.29\textwidth]{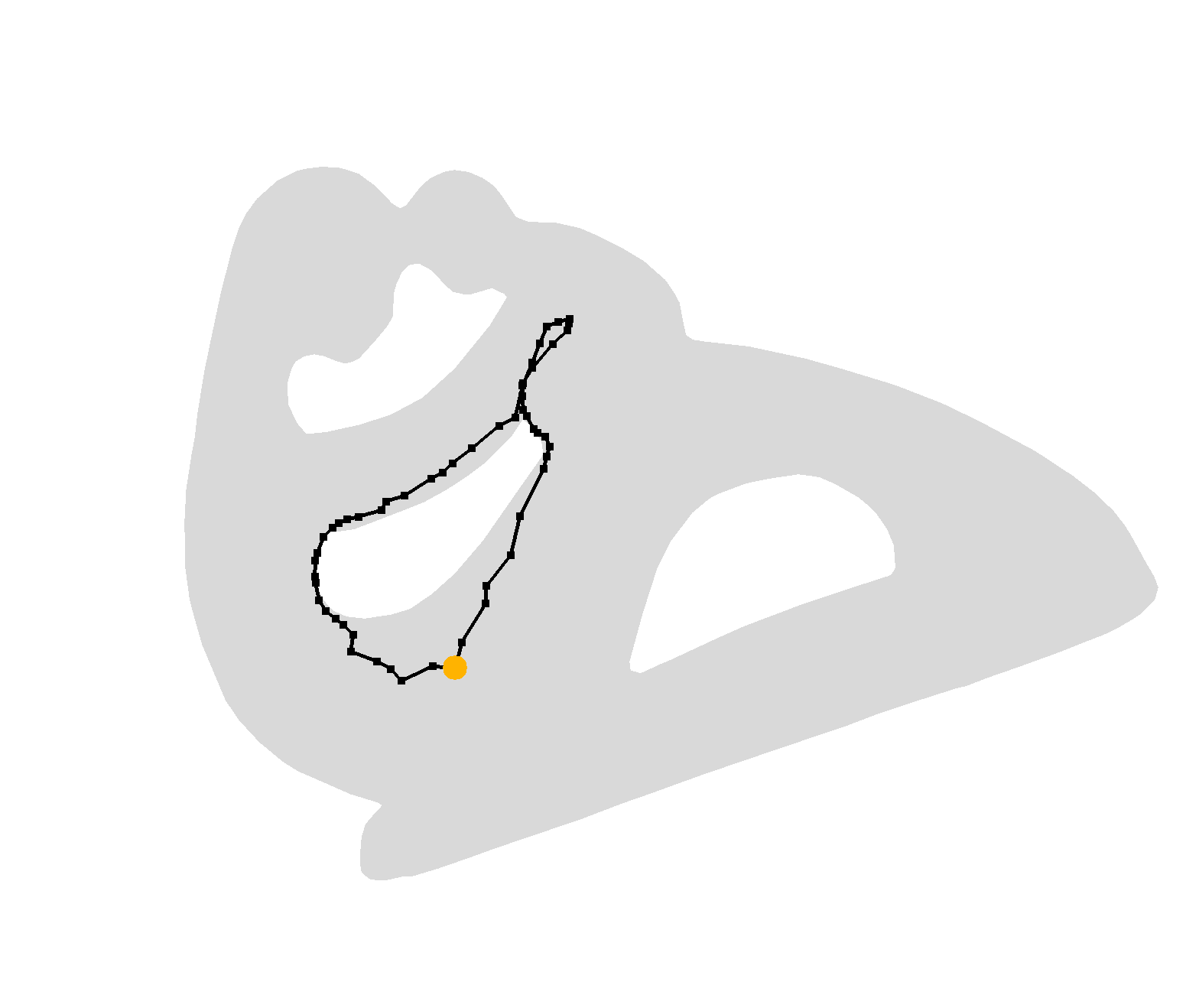}
    \includegraphics[width=0.29\textwidth]{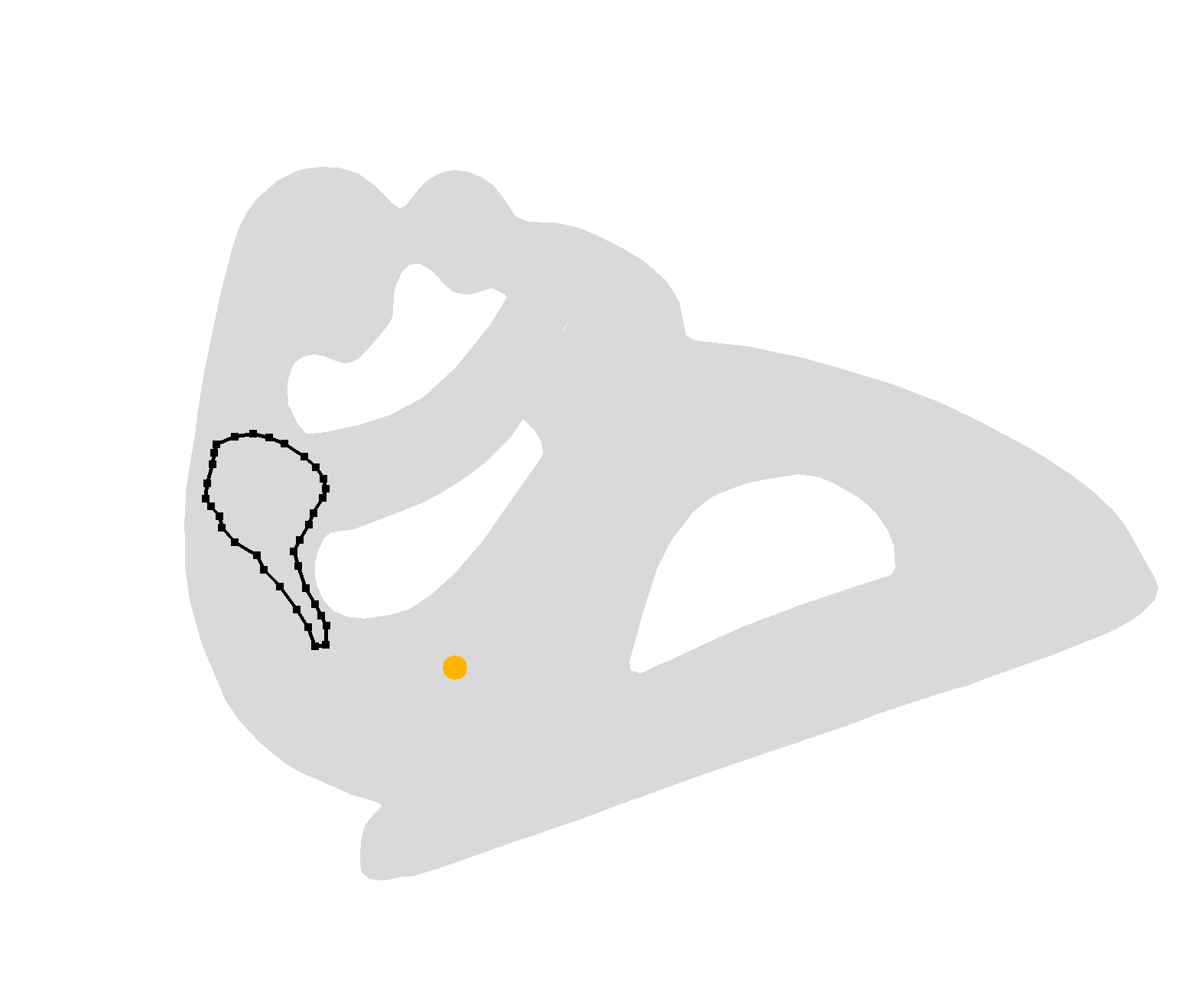}
    \includegraphics[width=0.29\textwidth]{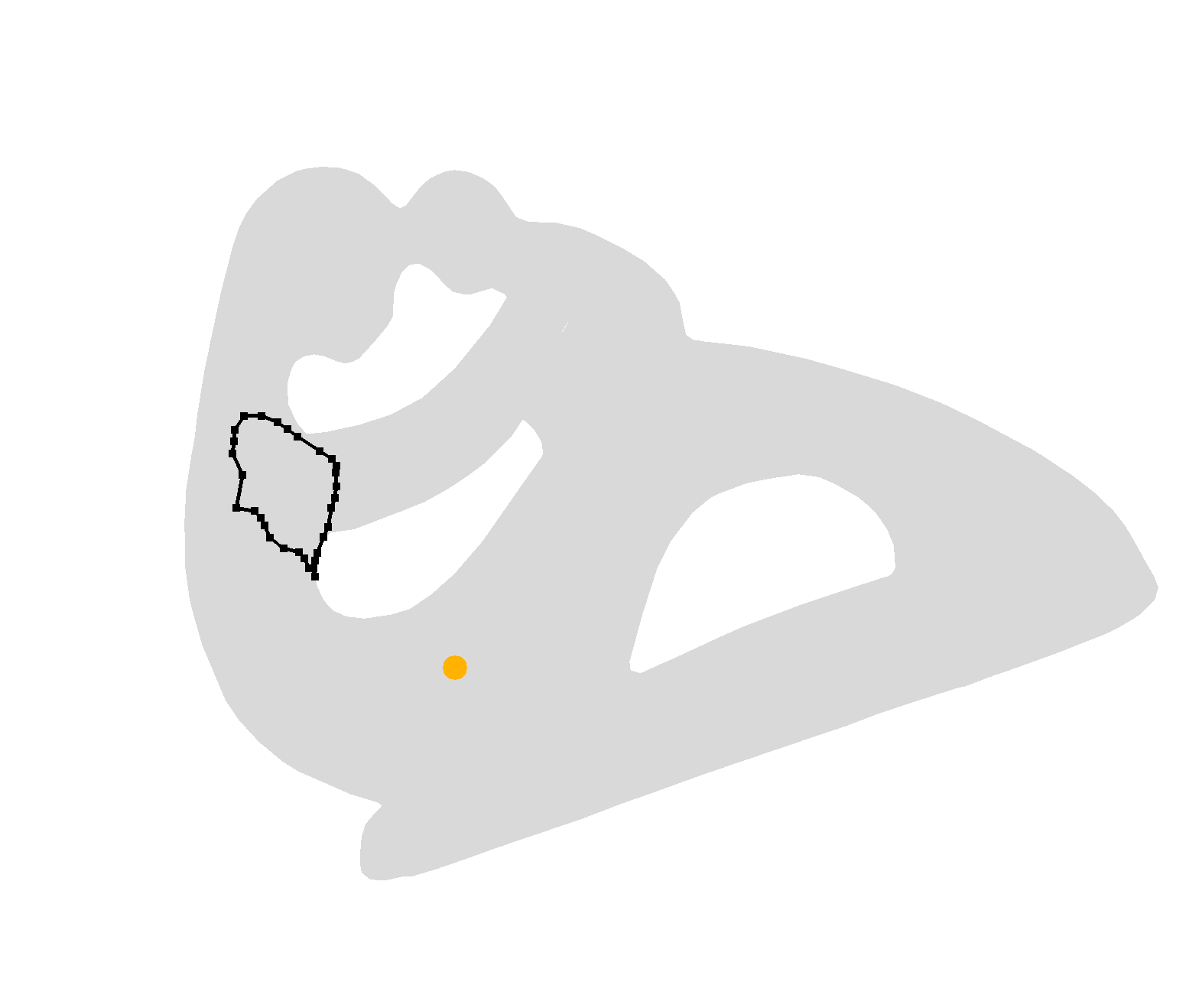}\\
    \caption{The fertility model and a selected vertex (top left). The remaining eight images show the generator curves (on top of a transparent rendering of the mesh) produced with this vertex as the starting point for the shortest-path tree.}
    \label{fig:fertility-generators}
\end{figure}

\begin{figure}[p]
    \centering
    \includegraphics[height=5cm]{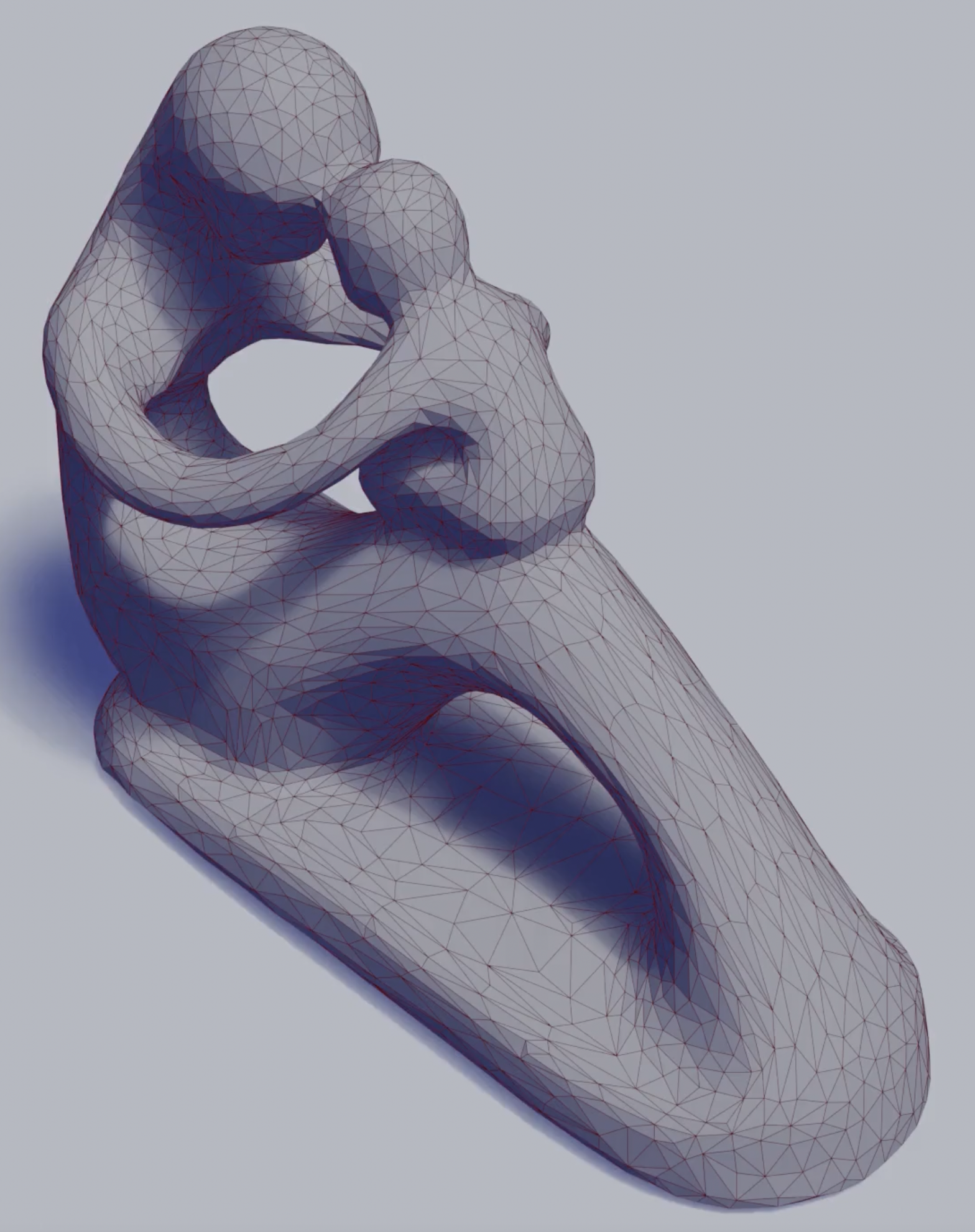}
    \includegraphics[height=5cm]{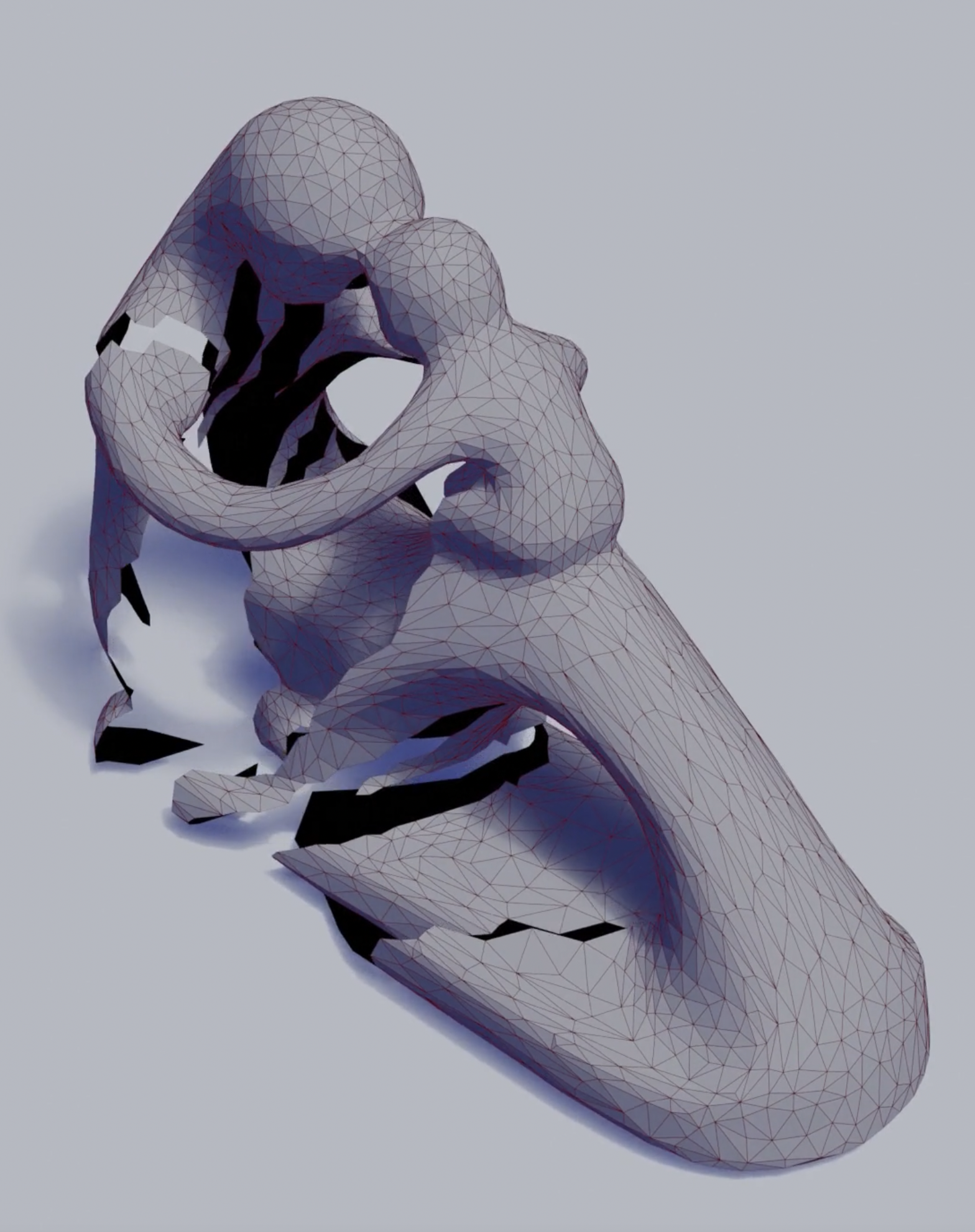}\\
    \includegraphics[height=5cm]{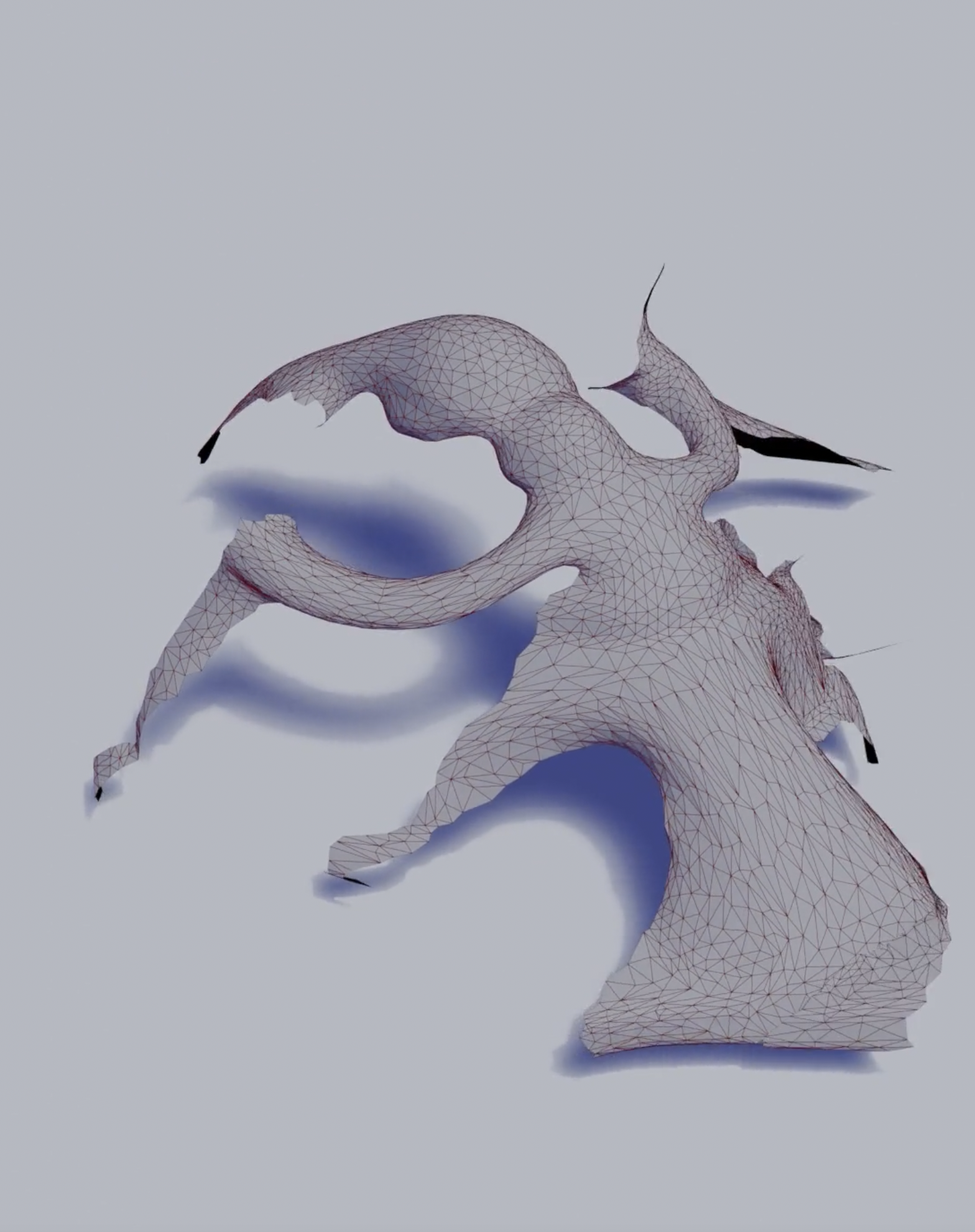}
    \includegraphics[height=5cm]{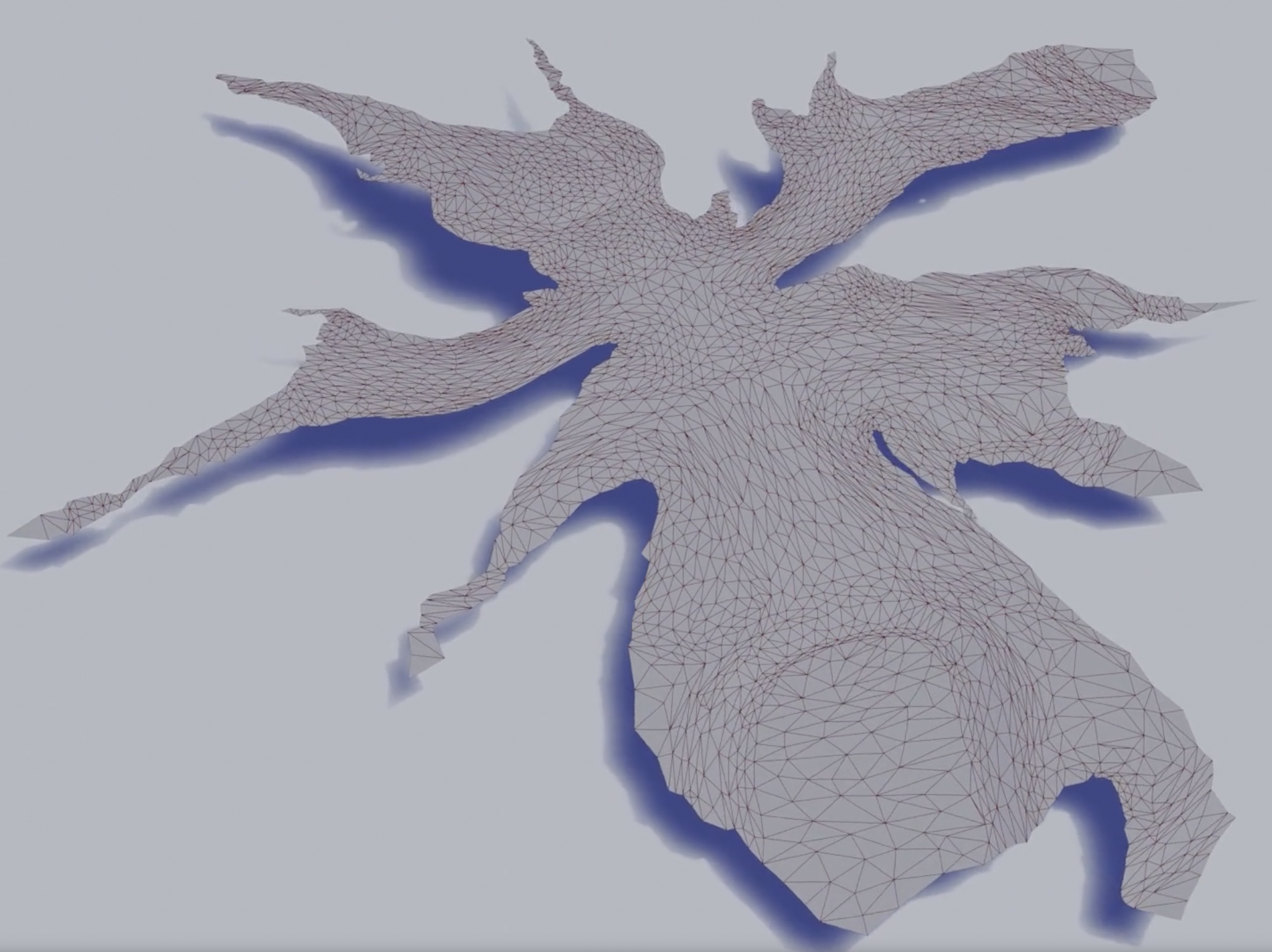}
    \caption{The fertility model is here cut along each of the eight generator curves and subsequently flattened.}
    \label{fig:fertility-cut}
\end{figure}

\marginnote{\footnotesize\itshape Note that the two figures are from different runs of the program, and the sets of generator curves are different.}Figure~\ref{fig:fertility-generators} shows the Fertility model together with a starting vertex (orange) and eight pictures that show the individual generator curves from the cut graph. In Figure~\ref{fig:fertility-cut} we see the Fertility model being flattened into a disc after it has been cut along the cut graph.

Note that this example was generated using the GEL library (\url{https://github.com/janba/GEL}). It is available in GEL as a Python example that can be run using PyGEL (\url{https://pypi.org/project/PyGEL3D/}).
\section{Conclusions}
In this document, we have described several topological properties, presented methods for computing the key quantitative properties -- $s$, $g$, and $b$ -- and derived the Betti numbers using just the Euler-Poincaré formula together with basic observations about the number of edges in a spanning tree. Finally, a method for computing the cut graph of a mesh was presented. Throughout the text, we have avoided relying on abstract algebra. The price we pay for this simplification is that the scope is narrow. The tools provided in this paper do not immediately generalise to structures other than polygonal meshes. Hopefully, what was lost in terms of generality was gained in terms of accessibility.
\section*{Acknowledgements}
Thanks to David Brander for helpful comments on nomenclature.
\bibliographystyle{unsrturl}
\bibliography{ref}
\end{document}